\documentclass[preprint,1p,10pt]{elsarticle}

\usepackage{graphicx}
\usepackage{amssymb}

\usepackage{amsmath}

\usepackage{color}

\newtheorem{remark}{Remark}

\biboptions{sort&compress,comma,round}

\journal{CMAME}

\begin{document}

\begin{frontmatter}


\title{Dispersion-optimized quadrature rules for isogeometric analysis: modified inner products, their dispersion properties, and optimally blended schemes}



\author[curtin]{Vladimir Puzyrev\corref{cor1}}
\ead{vladimir.puzyrev@gmail.com}

\author[curtin]{Quanling Deng}
\ead{qdeng12@gmail.com}

\author[curtin,csiro]{Victor Calo}
\ead{vmcalo@gmail.com}

\cortext[cor1]{Corresponding author}

\address[curtin]{Department of Applied Geology,
Western Australian School of Mines, Curtin University,
Kent Street, Bentley, Perth, WA 6102, Australia}
\address[csiro]{Mineral Resources,
Commonwealth Scientific and Industrial Research Organisation (CSIRO),
Kensington, Perth, WA 6152, Australia}

\begin{abstract}
This paper introduces optimally-blended quadrature rules for isogeometric analysis and analyzes the numerical dispersion of the resulting discretizations. To quantify the approximation errors when we modify the inner products, we generalize the Pythagorean eigenvalue theorem of Strang and Fix. The proposed blended quadrature rules have advantages over alternative integration rules for isogeometric analysis on uniform and non-uniform meshes as well as for different polynomial orders and continuity of the basis. The optimally-blended schemes improve the convergence rate of the method by two orders with respect to the fully-integrated Galerkin method. The proposed technique increases the accuracy and robustness of isogeometric analysis for wave propagation problems.
\end{abstract}


\end{frontmatter}


\section{Introduction}
\label{S:1}

Dispersion analysis is a powerful tool to understand the approximation errors of a numerical method. Amongst the most popular numerical methods used for wave propagation problems are the finite element and the spectral element methods \cite{komatitsch1998spectral, komatitsch1999introduction, de2010stability}, whose implementations exhibit excellent parallel scalability \cite{burstedde2010extreme}. The dispersive properties of these methods have been studied in detail to suggest that the most cost-effective scheme can be obtained by an appropriate weighted average of both methods \citep{marfurt1984accuracy, Ainsworth2010, seriani2007optimal}. This idea was employed in \citep{ainsworth2010optimally} where the optimal blending of the finite element and the spectral element methods was obtained and shown to provide two orders of extra accuracy (superconvergence) in the dispersion error. The blended scheme is equivalent to the use of nonstandard quadrature rules and therefore it can be efficiently implemented by replacing the standard Gaussian quadrature by a nonstandard rule \citep{ainsworth2010optimally, Ainsworth2010}.

This idea can be extended beyond the finite and spectral elements. The interest and development of isogeometric analysis (IGA) methods for partial differential equations has been continuously increasing in the last decade (see \citep{hughes2005isogeometric, bazilevs2007variational, bazilevs2010isogeometric, cottrell2006isogeometric, cottrell2007studies, hughes2008duality, hughes2010efficient, hughes2014finite, gomez2008isogeometric, bazilevs2010isogeometric2, gomez2010isogeometric, motlagh2013simulation, vignal2013phase, Vignal2016, calo2008simulation, zhang2007patient, bazilevs2009patient, bazilevs20113d, kamensky2015immersogeometric, liu2013isogeometric, thomas2015bezier, elguedj2008and, lipton2010robustness, elguedj2014isogeometric} and references therein). Isogeometric analysis uses as the basis functions those employed in computer aided design (CAD) systems that are capable of representing various complex geometries exactly. The initial motivation behind isogeometric analysis was to bridge the gap existing between CAD and the finite element analysis (FEA) framework, in particular to simplify the mesh generation process and the refinement of this mesh, as well as the transfer of the simulation information back to the design process. Additionally, isogeometric analysis has other attractive features. Its basis functions can have higher continuity across element interfaces. These functions may have up to $p-1$ continuous derivatives across element boundaries, where $p$ is the order of the underlying polynomial. Consequently, the approximated solutions have global continuity of order up to $p-1$. This local control of the continuity of the basis is a powerful tool in isogeometric analysis \citep{garcia2016value}. The publications \citep{cottrell2006isogeometric, cottrell2007studies, hughes2008duality, hughes2010efficient, hughes2014finite, reali2004isogeometric} show that highly continuous isogeometric analysis delivers more robustness and better accuracy per degree of freedom than standard finite elements, though at the cost of increased computational complexity \citep{collier2012cost, collier2013cost, collier2014computational, pardo2012survey}. To address the increased solution cost, several design strategies have been proposed for direct solvers \cite{wozniak2014computational, wozniak2015computational, los2015dynamics, kuznik2012graph}.

We analyze the errors in the discrete approximation to linear boundary- and initial-value problems by expressing them in terms of the eigenvalue and eigenfunction errors of the corresponding eigenproblem. This analysis requires a global error description, i.e. a characterization of the errors in the eigenvalues and the eigenfunctions for all the modes. Following \citet{hughes2014finite}, the total ``error budget'' determines the errors of the numerical method which consist of the eigenvalue and eigenfunction approximation errors. At each degree of freedom, the sum of the eigenvalue error and the square of the $L_2$-error in the eigenfunction scaled by the exact eigenvalue equals the square of the error in the energy norm. To quantify the approximation errors in the case when the integrals in the Galerkin formulation are not fully integrated, we generalize Strang's Pythagorean theorem to include the effect of inexact integration. We generalize the Pythagorean eigenvalue theorem for modified inner products.

Our study on blended quadratures seeks to increase the accuracy and robustness of the isogeometric framework. Smoother, such as $C^{p-1}$ continuous, basis functions produce better approximations of derivatives than $C^0$ finite elements in second-order problems \citep{hughes2008duality}. We focus on dispersion analysis and show that a modified inner product can be obtained by adapting the blending ideas introduced by \citet{ainsworth2010optimally} to this highly-continuous discretizations. Thus, the new blending schemes can be used in conjunction with smoother numerical methods to reduce the errors in the approximation of the eigenvalues (and, in some cases, the eigenfunctions). We find that using optimal blendings improves the convergence rate from the optimal for the polynomial order by two orders while not increasing the solution cost.

We consider an elliptic eigenvalue problem which describes the normal modes and frequencies of free structural vibrations, a classical application in engineering, which has been thoroughly studied in the isogeometric analysis literature. A similar dispersion analysis has been applied to the Helmholtz equation for time-harmonic wave propagation that arises in acoustics and electromagnetics.

The outline of this work is as follows. We first describe the model problem under study in Section 2. In Section 3, we present a generalization of the Pythagorean eigenvalue theorem that accounts for modified inner products. This generalization describes the terms contributing to the total error of a compatible numerical method. In Section 4, we describe the optimal blending of finite and spectral elements and present an optimal blending scheme for isogeometric analysis. Alternative blending strategies are discussed in Section 5. Numerical examples for one-dimensional and two-dimensional problems are given in Section 6. We perform numerical analysis of the discrete frequency spectra for the one- and two-dimensional cases and show the eigenvalue and eigenfunction errors in the $L_2$ and energy norms. Finally, Section 7 summarizes our findings and describes future research directions.

\section{Problem statement}
\label{S:2}

We start our analysis with the Laplace eigenvalue problem
\begin{align} \label{eq:Original}
\begin{split}
& \Delta u = \lambda u \ \ \textrm{in} \ \Omega \\
& u=0 \ \ \textrm{on} \ \partial \Omega,
\end{split}
\end{align}
where $\Delta = \nabla^2$ is the Laplace operator. In particular, we analyze the one-dimensional (1D)
\begin{align} \label{eq:Original1D}
\begin{split}
& \frac{\partial^2 u}{\partial x^2} = \lambda u,\ \ \textrm{for}\ x \in \Omega= ]0,1[ \\
& u(0) = u(1) = 0,
\end{split}
\end{align}
and the two-dimensional (2D)
\begin{align} \label{eq:Original2D}
\begin{split}
& \frac{\partial^2 u}{\partial x^2} + \frac{\partial^2 u}{\partial y^2} = \lambda u,\ \ \textrm{for}\ (x,y) \in \Omega = \ ]0,1[\ \times \ ]0,1[ \\ 
& u(x,y) \rvert _{\partial \Omega} = 0,
\end{split}
\end{align}
eigenvalue problems with homogeneous Dirichlet boundary conditions. The following derivations are specialized for the one-dimensional case to simplify the discussion. Such a problem has a countable infinite set of eigenvalues ${\lambda _j} \in {\rm{R}}$ and an associated set of orthonormal eigenfunctions ${u_j}$
\begin{equation}
{\lambda _1} \le {\lambda _2} \le ... \le {\lambda _j} \le ...
\end{equation}
\begin{equation}
({u_j},{u_k}) = \int_{\Omega}  {{u_j}(x){u_k}} (x)dx = {\delta _{jk}},
\end{equation}
where $\delta _{jk}$ is the Kronecker delta which is equal to 1 when $i=j$ and 0 otherwise. The eigenvalues are real, positive, and countable, for each eigenvalue $\lambda _j$ there exists an eigenfunction $u_j$. The normalized eigenfunctions form an $L_2$-orthonormal basis and they are orthogonal also in the energy inner product
\begin{equation}
(\mathcal{L}{u_j},{u_k}) = ({\lambda _j}{u_j},{u_k}) = {\lambda _j}{\delta _{jk}}.
\end{equation}

The standard weak form for the eigenvalue problem is stated as follows: Find all eigenvalues ${\lambda _j} \in {\rm{R}}$ and eigenfunctions ${u_j} \in V$ such that, for all $w \in V$,
\begin{equation} \label{eq:weak}
a(w,{u_j}) = {\lambda _j}(w,{u_j}),
\end{equation}
where
\begin{equation}
a(w,{u_j}) = \int_{\Omega}  {\frac{d w}{dx}\frac{d u_j}{dx}} dx,
\end{equation}
and $V$ is a closed subspace of ${{H^1}(\Omega )}$. We use the notation of \citet{strang1973analysis} where $( \cdot , \cdot )$ and $a( \cdot , \cdot )$ are symmetric bilinear forms which define the following inner products
\begin{equation}
\left\| w \right\|_E^2 = a(w,w),
\end{equation}
\begin{equation}
\left\| w \right\|_{}^2 = (w,w),
\end{equation}
for all $v,w \in V$. The energy norm is denoted as ${\left\| {\; \cdot \;} \right\|_E}$ and is equivalent to the ${{H^1}(\Omega )}$ norm on $V$ and $\left\| {\; \cdot \;} \right\|$ is the standard ${{{L_2}(\Omega )} = {H^0}(\Omega ) }$ norm. Here ${{H^1}(\Omega )}$ denotes the Sobolev space of functions
\begin{equation}
{{H^1}(\Omega )} = \{ f : \Omega \rightarrow R \ \vert \  \left\| {f} \right\| _{H^1} < \infty \}, \ \ \ \left\| {f} \right\| ^{2} _{H^1} = \int\limits_a^b \left[ f^2(x) + \left( \frac{d}{dx} f(x) \right)^2 \right] dx.
\end{equation}

The Galerkin formulation of the eigenvalue problem is the discrete form of \eqref{eq:weak}: Find the discrete eigenvalues $\lambda _j^h \in {\rm{R}}$ and eigenfunctions $u_j^h \in V^h \subset V$ such that, for all ${w^h} \in {V^h} \subset V$,
\begin{equation} \label{eq:galerkin}
a({w^h},u_j^h) = \lambda _j^h({w^h},u_j^h).
\end{equation}

Because of the minimax principle \citep{strang1973analysis}, all discrete eigenvalues resulting from the finite element method approximate from above the exact eigenvalues
\begin{equation} \label{eq:above}
\lambda _j^h \ge \lambda _j \text{\ \ \ \ for all} \ j .
\end{equation}

The theorem for the Pythagorean eigenvalue error can be applied in this case. For each discrete mode, the eigenvalue error and the product of the eigenvalue and the square of the eigenfunction error in the $L_2$-norm sum to the square of the error in the energy norm
\begin{equation} \label{eq:ptorig1}
{\lambda _j^h - {\lambda _j}} + {\lambda _j}{\left\| {u_j^h - {u_j}} \right\|}^2 = {\left\| {u_j^h - {u_j}} \right\|_E^2},
\end{equation}
which can be rewritten as
\begin{equation} \label{eq:ptorig12}
\frac{{\lambda _j^h - {\lambda _j}}}{{{\lambda _j}}} + \frac{{{{\left\| {u_j^h - {u_j}} \right\|}^2}}}{{{{\left\| {{u_j}} \right\|}^2}}} = \frac{{\left\| {u_j^h - {u_j}} \right\|_E^2}}{{\left\| {{u_j}} \right\|_E^2}}.
\end{equation}

By using the standard normalization \citep{strang1973analysis, hughes2014finite}, that is $\left\| {{u_j}} \right\| = 1$ and using \eqref{eq:weak}, \eqref{eq:ptorig12} simplifies to
\begin{equation} \label{eq:ptorig2}
\frac{{\lambda _j^h - {\lambda _j}}}{{{\lambda _j}}} + {\left\| {u_j^h - {u_j}} \right\|^2} = \frac{{\left\| {u_j^h - {u_j}} \right\|_E^2}}{{{\lambda _j}}}.
\end{equation}

Since the second term is always non-negative and taking into account \eqref{eq:above}, the following inequalities can be deduced from \eqref{eq:ptorig2}
\begin{equation}
\lambda _j^h - {\lambda _j} \le \left\| {u_j^h - {u_j}} \right\|_E^2
\end{equation}
\begin{equation} \label{eq:ineq1}
{\left\| {u_j^h - {u_j}} \right\|^2} \le \frac{{\left\| {u_j^h - {u_j}} \right\|_E^2}}{{{\lambda _j}}}.
\end{equation}

The inequality \eqref{eq:ineq1} does not hold for methods that do not approximate all eigenvalues from above, for example, the spectral element method \citep{ainsworth2009dispersive}. If the discrete method does not fully reproduce the inner product, as the theorem assumes, the theorem needs to be adapted to account for the error in the inner product representation at the discrete level.

The original eigenvalue problem \eqref{eq:Original} can be written in matrix form as
\begin{equation}
\mathbf{K} \boldsymbol{\phi_j} = {\lambda_j} \mathbf{M} \boldsymbol{\phi_j}
\end{equation}
and the global mass and stiffness matrices $\mathbf{M}$ and $\mathbf{K}$ are assembled from the elemental matrices that are given by
\begin{equation} \label{eq:mass}
\mathbf{M}_{i j} = \int\limits_{a}^{b} \phi_i(x)  \phi_j(x) dx
\end{equation}
\begin{equation} \label{eq:stiff}
\mathbf{K}_{i j} = \int\limits_{a}^{b} \frac{d \phi_i}{dx}  \frac{d \phi_j}{dx} dx ,
\end{equation}
where $\phi_i(x)$ are the piecewise polynomial basis functions. $\mathbf{M}$ and $\mathbf{K}$ are symmetric positive definite matrices and in the 1D case they are banded with the maximum bandwidth equal to $p$.

We give a brief description to introduce notation and to analyze the effects of quadrature for finite elements, spectral elements, and isogeometric analysis. The integrals in the elemental matrices \eqref{eq:mass} and \eqref{eq:stiff} are evaluated with the help of numerical integration rules. An $m$-point quadrature rule requires $m$ evaluations of a function $f(x)$ to approximate its weighted integral over an interval $[a, b]$
\begin{equation}
\int_{a}^{b} w(x) f(x) dx = \sum_{i=1}^{m} w_i f(x_i) + R_m.
\end{equation}
Here $w_i$ and $x_i$ are the weights and the nodes of the rule, respectively. The quadrature rule is exact for a given function $f(x)$ when the remainder $R_m$ is exactly zero. For example, the standard $m$-point Gauss-Legendre (GL or Gauss) quadrature is exact for the linear space of polynomials of degree at most $2m-1$ (see, for example, \cite{stoer2013introduction, bartovn2016optimal}).

The classical Galerkin finite element analysis typically employs the Gauss quadrature with $p+1$ quadrature points per parametric direction that fully integrates every term in the bilinear forms defined by the weak form. A quadrature rule is optimal if the function is evaluated with the minimal number of nodes (for example, Gauss quadrature with $m$ evaluations is optimal for polynomials of order $2m-1$ in one dimension).

Element-level integrals may be approximated using other quadrature rules, for example the Gauss-Lobatto-Legendre (GLL or Lobatto) quadrature rule that is used in the spectral element method (SEM). The Lobatto quadrature evaluated at $m$ nodes is accurate for polynomials up to degree $2m-3$. However, selecting a rule with $p+1$ evaluations for a polynomial of order $p$ and collocating the Lagrange nodes with the quadrature positions renders the mass matrix diagonal in 1D, 2D and 3D for arbitrary geometrical mappings. This resulting diagonal mass matrix is a more relevant result than the reduction in the accuracy of the calculation. Particularly, given that this property preserves the optimal convergence order for these higher-order schemes. Lastly, the spectral elements possess a superior phase accuracy when compared with the standard finite elements of the same polynomial order \citep{ainsworth2009dispersive}.

Isogeometric analysis based on NURBS (Non-Uniform Rational B-Splines) has been described in a number of papers (e.g. \cite{cottrell2006isogeometric, bazilevs2007variational, cottrell2007studies, hughes2008duality}) and the efficient implementation of the method in open source software has been discussed in \citep{de2011geopdes, pauletti2015igatools, dalcin2016petiga, sarmiento2016petiga}. Isogeometric analysis employs piecewise polynomial curves composed of linear combinations of B-spline basis functions. B-spline curves of polynomial order $p$ may have up to $p-1$ continuous derivatives across element boundaries. Three different refinement mechanisms are commonly used in isogeometric analysis, namely the \textit{h}-, \textit{p}- and \textit{k}-refinement, as detailed in \cite{cottrell2007studies}. We refer the reader to \citet{piegl2012nurbs} for the definition of common concepts of isogeometric analysis such as knot vectors, B-spline functions, and NURBS.

The derivation of optimal quadrature rules for NURBS-based isogeometric analysis with spaces of high polynomial degree and high continuity has attracted significant attention in recent years \citep{hughes2010efficient, auricchio2012simple, ait2015explicit, bartovn2015gaussian, bartovn2016optimal, bartovn2016gaussian, calabro2016fast, Adam2015732, antolin2015efficient, hiemstra2016optimal}. The efficiency of Galerkin-type numerical methods for partial differential equations depends on the formation and assembly procedures, which, in turn, largely depend on the efficiency of the quadrature rule employed. Integral evaluations based on full Gauss quadrature are known to be efficient for standard $C^0$ finite element methods, but inefficient for isogeometric analysis that uses higher-order continuous spline basis functions \citep{nguyen2016collocated}.

\citet{hughes2010efficient} studied the effect of reduced Gauss integration on the finite element and isogeometric analysis eigenvalue problems. By using $p$ Gauss points (i.e., underintegrating using one point less), one modifies the mass matrix only (in 1D). By using less than $p$ Gauss points (i.e., underintegrating using several points less), both mass and stiffness matrices are underintegrated. Large underintegration errors may lead to the loss of stability since the stiffness matrix becomes singular. As shown in \cite{hughes2010efficient}, this kind of underintegration led to the results that were worse than the fully integrated ones and the highest frequency errors diverged as the mesh was refined. However, as we show in the next sections, using properly designed alternative quadratures may lead to more accurate results.

The assembly of the elemental matrices into the global stiffness and mass matrices is done in a similar way for all Galerkin methods we analyze in this paper. Similarly, the convergence rate for all Galerkin schemes we analyze is the same. However, the heterogeneity of the high-order finite element ($C^0$ elements, i.e., SEM and FEA) basis functions leads to a branching of the discrete spectrum and a fast degradation of the accuracy for higher frequencies. In fact, the degraded frequencies in 1D are about half of all frequencies, while in 3D this proportion reduces to about seven eighths. On uniform meshes, B-spline basis functions of the highest $p-1$ continuity, on the contrary, are homogeneous and do not exhibit such branching patterns other than the outliers that correspond to the basis functions with support on the boundaries of the domain.

\begin{remark}
The above analysis also applies to the Helmholtz equation (e.g. \cite{hughes2008duality}). The classical wave propagation equation is
\begin{equation}
\Delta u - \frac{1}{c^2} \frac{\partial^2 u}{\partial^2 t} = 0 .
\end{equation}
Assuming time-harmonic solutions of the form $u(\mathbf{x}, t) = e^{−iwt}u(\mathbf{x})$ for a given temporal frequency $\omega$, the wave equation reduces to the Helmholtz equation

\begin{equation} \label{eq:Helm}
\Delta u + k^2 u = 0 ,
\end{equation}
where the wavenumber $k=\omega/c$ represents the ratio of the angular frequency $\omega$ to the wave propagation speed $c$. The wavelength is equal to $2\pi/k$.

The discretization of \eqref{eq:Helm} leads to the following linear equation system
\begin{equation} \label{eq:discrete}
\left( \mathbf{K} - k^2 \mathbf{M} \right) \boldsymbol{u}^h = 0,
\end{equation}
and the relation between the eigenvalues and wavenumbers is $\lambda_j = k_j^2$. 

The solution of \eqref{eq:discrete} is a linear combination of plane waves with numerical wavenumbers $k^h$, where $k^h \neq k$ and hence the discrete and exact waves have different wavelengths. The goal of the dispersion analysis is to quantify this difference and define this difference as the dispersion of the numerical method, i.e. how well the discrete wavenumber $k^h$ approximates the continuous $k$ \citep{strang1971finite, babuvska1989finite, hughes2008duality}.
\end{remark}

\begin{remark}
For multidimensional problems on tensor product grids, the stiffness and mass matrices can be expressed as Kronecker products of 1D matrices \cite{gao2014fast, gao2015preconditioners}. For example, in the 2D case, the components of $\mathbf{K}$ and $\mathbf{M}$ can be represented as fourth-order tensors using the definitions of the matrices and the basis functions for the 1D case \citep{de2007grid, gao2014fast}
\begin{equation}
\mathbf{M} _{i j k l} = \mathbf{M} ^{1D} _{i k} \mathbf{M} ^{1D} _{j l}
\end{equation}
\begin{equation}
\mathbf{K} _{i j k l} = \mathbf{K} ^{1D} _{i k} \mathbf{M} ^{1D} _{j l} + \mathbf{K} ^{1D} _{j l} \mathbf{M} ^{1D} _{i k},
\end{equation}
where $\mathbf{M} ^{1D} _{i j}$ and $\mathbf{K} ^{1D} _{i j}$ are the mass and stiffness matrices of the 1D problem as given by \eqref{eq:mass}, \eqref{eq:stiff}. We refer the reader to \citep{de2007grid} for the description of the summation rules.
\end{remark}

\section{Pythagorean eigenvalue theorem for weak forms with modified inner products}
\label{S:3}

To understand the accuracy delivered by each quadrature technique, we measure the error they induce in the inner product. In this section we generalize the Pythagorean eigenvalue error theorem to account for modified inner products. From \eqref{eq:weak} and \eqref{eq:galerkin}, the exact and fully-integrated approximate solutions are 
\begin{equation}
a({u_j},{u_j}) = \lambda _j ({u_j},{u_j})
\end{equation}
\begin{equation}
a(u_j^h,u_j^h) = \lambda _j^h(u_j^h,u_j^h).
\end{equation}

We write the approximate solutions for the modified inner-product discretizations as
\begin{equation}
{a_h}(v_j^h,v_j^h) = \mu _j^h{(v_j^h,v_j^h)_h},
\end{equation}
where $a_h( \cdot , \cdot )$ and $( \cdot , \cdot )_h $ are the discrete inner products represented by the chosen quadratures. For each $j$, the set $v_j^h$ and $\mu_j^h$ are the discrete eigenpair resulting from these discrete inner product representations.

In all these three cases, the corresponding eigenfunctions are orthonormal with respect to the appropriate inner product
\begin{equation}
({u_i},{u_j}) = (u_i^h,u_j^h) = {(v_i^h,v_j^h)_h} = {\delta _{ij}}.
\end{equation}

The Pythagorean eigenvalue theorem for the fully integrated case where the discrete inner product is equivalent to the continuous one for all functions in the finite dimensional spaces used, is derived as follows \cite{strang1973analysis}
\begin{align} \label{eq:pt0}
\left\| {{u_j} - u_j^h} \right\|_E^2 &= a({u_j} - u_j^h,{u_j} - u_j^h) \nonumber\\
&= a({u_j},{u_j}) - 2a({u_j},u_j^h) + a(u_j^h,u_j^h) \nonumber\\
&= {\lambda _j} - 2{\lambda _j}({u_j},u_j^h) + \lambda _j^h \nonumber\\
&= {\lambda _j}\left[ {1 - 2({u_j},u_j^h) +1} \right] + \lambda _j^h - {\lambda _j} \nonumber\\
&= {\lambda _j}\left[ {({u_j},{u_j}) - 2({u_j},u_j^h) + (u_j^h,u_j^h)} \right] + \lambda _j^h - {\lambda _j} \nonumber\\
&= {\lambda _j}\left\| {{u_j} - u_j^h} \right\|_0^2 + \lambda _j^h - {\lambda _j},
\end{align}
where the first equality is valid due to the definition of $a( \cdot , \cdot )$ as the energy norm, the second one is due to bilinearity of $a( \cdot , \cdot )$, the third one is due to the definition of the continuous eigenproblem, the fourth one is a simple regrouping of terms, the fifth is due to the $L_2$ norms scaling of the continuous and discrete eigenfunctions, while the sixth equality is another regrouping of terms to yield the final result.

In the modified inner product case we can write the generalized Pythagorean theorem as follows
\begin{align} \label{eq:pt1}
\begin{split}
\left\| {{u_j} - v_j^h} \right\|_E^2 = {\lambda _j}\left[ {1 - 2({u_j},v_j^h)} +1 \right] + a(v_j^h,v_j^h) - {\lambda _j},
\end{split}
\end{align}
where the first four steps of \eqref{eq:pt0} apply unchanged. Since $\mu _j^h{(v_j^h,v_j^h)_h} - {a_h}(v_j^h,v_j^h) = 0$, we can add this to \eqref{eq:pt1} which results in
\begin{align}
\left\| {{u_j} - v_j^h} \right\|_E^2 =& \ {\lambda _j}\left[ {1 - 2({u_j},v_j^h)} +1 \right] + \mu _j^h{(v_j^h,v_j^h)_h} - {\lambda _j} + a(v_j^h,v_j^h) - {a_h}(v_j^h,v_j^h).
\end{align}

Then, using the orthonormality of the eigenfunctions in their respective inner products $({u_j},{u_j}) = 1$ and ${(v_j^h,v_j^h)_h} = 1$, while, in general, the following is true: $(v_j^h,v_j^h) \neq 1$. That is, it is not necessarily normalized to one in the continuous inner product, we finally obtain the modified theorem
\begin{align} \label{eq:pt2}
\begin{split}
\left\| {{u_j} - v_j^h} \right\|_E^2 &= \mu _j^h - {\lambda _j} + (a - {a_h})(v_j^h,v_j^h) \\
& + {\lambda _j}\left[ {({u_j},{u_j}) - 2({u_j},v_j^h) + (v_j^h,v_j^h)} \right] + {\lambda _j}\left[ {1 - (v_j^h,v_j^h)} \right] \\
& = \mu _j^h - {\lambda _j} + {\lambda _j}\left\| {{u_j} - v_j^h} \right\|_0^2 \\
&+ (a - {a_h})(v_j^h,v_j^h) + {\lambda _j}\left[ {{{(v_j^h,v_j^h)}_h} - (v_j^h,v_j^h)} \right].
\end{split}
\end{align}

The last equation in \eqref{eq:pt2} shows that the energy norm error for weak forms with modified inner products consists of four terms. The first two are similar to the Pythagorean eigenvalue theorem of \citet{strang1973analysis} and they represent the eigenvalue error and $L_2$-norm of the eigenfunction error. The third term of \eqref{eq:pt2} is

\begin{equation} \label{eq:term3}
(a - {a_h})(v_j^h,v_j^h) = \left\| {v_j^h} \right\|_E^2 - \left\| {v_j^h} \right\|_{E,h}^2.
\end{equation}
While the last term is

\begin{equation} \label{eq:term4}
{\lambda _j}\left[ {{{(v_j^h,v_j^h)}_h} - (v_j^h,v_j^h)} \right] = {\lambda _j}\left\| {v_j^h} \right\|_{0,h}^2 - {\lambda _j}\left\| {v_j^h} \right\|_0^2 = {\lambda _j} \left(1 - \left\| {v_j^h} \right\|_0^2 \right).
\end{equation}

The third term represents the error in the discrete energy norm, while the last one is the error in the $L_2$ inner product. In the numerical examples section, we illustrate the contribution of the two new terms to the total error of the numerical methods. In the cases when these terms are equal to zero, \eqref{eq:pt2} reduces to the standard Pythagorean eigenvalue error theorem which when the results of \eqref{eq:term3} and \eqref{eq:term4} are substituted in \eqref{eq:pt2} yield:
\begin{equation}
\left\| {{u_j} - v_j^h} \right\|_E^2 = \mu _j^h - {\lambda _j} + {\lambda _j}\left\| {{u_j} - v_j^h} \right\|_0^2 + \left\| {v_j^h} \right\|_E^2 - \left\| {v_j^h} \right\|_{E,h}^2 + {\lambda _j} \left(1 - \left\| {v_j^h} \right\|_0^2 \right),
\end{equation}
which can be simplified to
\begin{equation} \label{eq:ptfinal}
\frac{\left\| {{u_j} - v_j^h} \right\|_E^2}{\lambda _j} = \frac{\mu _j^h - {\lambda _j}}{\lambda _j} + \left\| {{u_j} - v_j^h} \right\|_0^2 + \frac{\left\| {v_j^h} \right\|_E^2 - \left\| {v_j^h} \right\|_{E,h}^2}{\lambda _j} + \left(1 - \left\| {v_j^h} \right\|_0^2 \right).
\end{equation}

\section{Optimal blending for finite elements and isogeometic analysis}
\label{S:4}

Several authors (e.g. \cite{seriani2007optimal, ainsworth2010optimally, Esterhazy2014830}) studied the blended spectral-finite element method that uses nonstandard quadrature rules to achieve an improvement of two orders of accuracy compared with the fully integrated schemes. This method is based on blending the full Gauss quadrature, which exactly integrates the bilinear forms to produce the mass and stiffness matrices, with the Lobatto quadrature, which underintegrates them. This methodology exploits the fact that the fully integrated finite elements exhibit phase lead when compared with the exact solutions, while the underintegrated with Lobatto quadrature methods, such as, spectral elements have phase lag.

\citet{ainsworth2010optimally} chose the blending parameter to maximize the order of accuracy in the phase error. They showed that the optimal choice for the blending parameter is given by weighting the spectral element and the finite element methods in the ratio $p$ to one scaled by $(p+1)$. As mentioned above, this optimally blended scheme improves by two orders the convergence rate of the blended method when compared against the finite or spectral element methods that were the ingredients used in the blending. The blended scheme can be realized in practice without assembly of the mass matrices for both schemes, but instead by replacing the standard Gaussian quadrature rule by an alternative rule, as Ainsworth and Wajid clearly explained in \cite{ainsworth2010optimally}. Thus, no additional computational cost is required by the blended scheme although the ability to generate a diagonal mass matrix by the underintegrated spectral method is lost.

To show how an improvement in the convergence rate is achieved, consider, for example, the approximate eigenfrequencies written as a series in $\Omega = \omega h$ for the linear finite and spectral elements, respectively \citep{ainsworth2010optimally}
\begin{equation} \label{eq:fe1}
\omega^h _ {FE} h = \Omega - \frac{\Omega^3}{24} + O(\Omega^5),
\end{equation}
\begin{equation} \label{eq:se1}
\omega^h _ {SE} h = \Omega + \frac{\Omega^3}{24} + O(\Omega^5).
\end{equation}

When these two schemes are blended using a blending parameter $\tau$, the approximate eigenfrequencies become
\begin{equation} \label{eq:fese1}
\omega^h _ {BL} h = \Omega + \frac{\Omega^3}{24}(2 \tau - 1) + O(\Omega^5).
\end{equation}

For $\tau = 0$ and $\tau = 1$, the above expression reduces to the ones obtained by the finite element and spectral element schemes, respectively. The choice of $\tau = 1/2$ allows the last term of \eqref{eq:fese1} to vanish and increases by two additional orders of accuracy the phase approximation when compared with the standard schemes. Similarly, by making the optimal choice of blending parameter $\tau = \frac{p}{(p+1)}$ in high-order schemes, they removed the leading order term from the error expansion.

The numerical examples in Section 6 show that a similar blending can be applied to the isogeometric mass and stiffness matrices to reduce the eigenvalue error. For $C^1$ quadratic elements, the approximate eigenfrequencies are
\begin{equation} \label{eq:iga21}
\omega^h _ {GL} h = \Omega - \frac{1}{5!}\frac{\Omega^5}{12} + O(\Omega^7),
\end{equation}
\begin{equation} \label{eq:iga22}
\omega^h _ {GLL} h = \Omega + \frac{1}{5!}\frac{\Omega^5}{24} + O(\Omega^7),
\end{equation}
and thus the optimal ratio of the Lobatto and Gauss quadratures is 2 : 1 ($\tau = 2/3$) similar to the optimally blended spectral-finite element scheme. For $C^2$ cubic elements, we determine that a non-convex blending with $\tau = 5/2$ allows us to remove the leading error term and thus achieve two additional orders of accuracy. For more details, we refer the reader to \cite{deng2017quadratures}.

(\ref{eq:fe1}-\ref{eq:iga22}) show that the \textit{absolute} errors in the eigenfrequencies converge with the rates of $O\left( \Omega^{2p+1} \right)$ and $O\left( \Omega^{2p+3} \right) $ for the standard and optimal schemes, respectively. If we consider the \textit{relative} eigenfrequency errors, from equations \eqref{eq:iga21} and \eqref{eq:iga22}, these take the form
\begin{equation} \label{eq:iga21rel}
\frac{\omega^h h}{\Omega} = 1 \pm \frac{\Omega^4}{\alpha} + ... ,
\end{equation}
that is, the convergence rate for frequencies computed using IGA approximations is $O\left( \Omega^{2p} \right) $ as shown in \cite{cottrell2006isogeometric, reali2004isogeometric}. The optimal blending in IGA leads to a $O\left( \Omega^{2p+2} \right) $ convergence rate for the relative eigenfrequencies. This superconvergence gain is similar to the one achieved by the optimally-blended spectral element method of \citet{ainsworth2010optimally}.

\section{Alternative blending schemes}
\label{S:5}

The optimal blending of \cite{ainsworth2010optimally} provides the best eigenvalue approximation when $\Omega \rightarrow 0$. However, in practical applications one might need to minimize the errors for different values of the wavenumber (frequency) for a set grid size rather than for very small ones (e.g., for a fixed value of grid points per wavelength in wave propagation problems) where the accuracy is already sufficiently high. \citet{seriani2007optimal} studied the possibility of blending the consistent and mass-lumped operators in spectral element schemes using a criteria that the phase error vanishes at a particular, user-specified, value of the normalized wavenumber. Certainly, in such a case the blending parameter is wavenumber and mesh dependent.

In the following numerical examples, we show that the alternative blendings in IGA lead to much smaller errors in the middle part of spectra. Our scheme blends the Gauss and Lobatto quadratures and we determine the blending parameter $\tau$ from numerical simulations to minimize an averaged error for a range of frequencies. In this context, we can also choose $\tau$ to minimize an expansion similar to \eqref{eq:fese1} for a given value of $\Omega$ if all relevant higher-order terms are available.

These nonstandard blendings lead to better approximations in the middle part of the spectrum for all orders of optimally blended IGA schemes. Finding explicit expressions of the equivalent quadrature rules is an important open problem which when solved will allow for efficient computation of the integration without resorting to computing with two different quadratures.

\section{Numerical examples}
\label{S:6}

In our numerical tests, we consider the one- and two-dimensional problems described in Section 2. The mesh is uniform unless otherwise specified and is chosen in such a way that $N$, the total number of degrees of freedom (or discrete modes), is $1000$ and $10000 \ (=100^2)$ in the 1D and 2D plots, respectively. 

\subsection{Optimally-blended methods in 1D}

The exact eigenvalues and corresponding eigenfunctions of the 1D problem are 
\begin{equation}
{\lambda _j} = {j^2}{\pi ^2},\ \ \ {u_j} = \sqrt 2 \sin (j\pi x),
\end{equation}
for $j = {1, 2, ... }$. The approximate eigenvalues $\lambda _j^h$ are sorted in ascending order and are compared to the corresponding eigenvalues of the exact operator ${\lambda _j}$.

In the following figures, we present the relative eigenvalue (EV) errors $\frac{\mu _j^h - {\lambda _j}}{\lambda _j}$, the $L_2$-norm eigenfunction (EF) errors $\left\| {{u_j} - v_j^h} \right\|_0^2$ and the relative energy-norm EF errors $\frac{\left\| {{u_j} - v_j^h} \right\|_E^2}{\lambda _j}$. This format of error representation clearly illustrates the budget of the generalized Pythagorean eigenvalue theorem  \eqref{eq:ptfinal}. The error in the $L_2$ norm $1 - \left\| {v_j^h} \right\|_0^2$ and the relative discrete energy norm error $\frac{\left\| {v_j^h} \right\|_E^2 - \left\| {v_j^h} \right\|_{E,h}^2}{\lambda _j}$ are shown whenever they are not zero.

\begin{figure}[!ht]
\centering\includegraphics[width=1.0\linewidth]{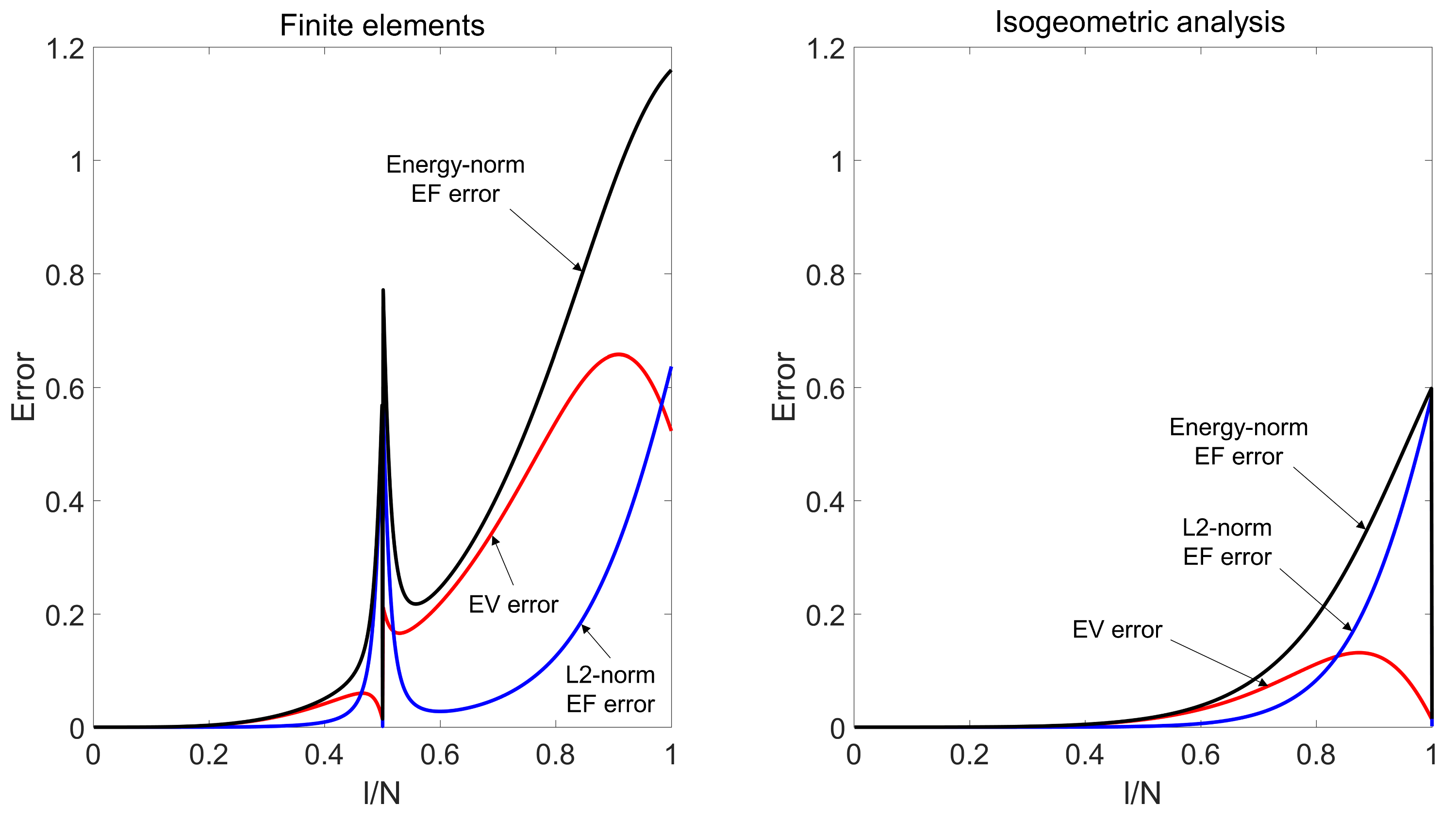}
\caption{Approximation errors for quadratic $C^0$ finite elements (left) and $C^1$ isogeometric elements (right), where EV and EF stand for eigenvalue and eigenfunction, respectively.}
\end{figure}

Figure 1 shows the approximation errors for quadratic finite elements ($C^0$) and isogeometric elements ($C^1$) using the standard Gaussian quadratures. Similar plots have been previously shown by \citet{cottrell2006isogeometric} and \citet{hughes2008duality, hughes2014finite} among others to illustrate that the approximate eigenvalues and eigenfunctions are significantly more accurate for IGA than for FEA for similar spatial resolutions. This improvement in the spectral accuracy of isogeometric analysis grows even larger for higher-order approximations. The large spikes in the eigenfunction errors that appear at the transition points between the acoustic and optical branches of the $C^0$ spectra, are absent in the maximum continuity discretizations. The figure shows that the standard Pythagorean eigenvalue error theorem is valid in this case. The last two terms of the modified Theorem \eqref{eq:ptfinal} are zero for all modes and not shown in this figure.

\begin{figure}[!ht]
\centering\includegraphics[width=1.0\linewidth]{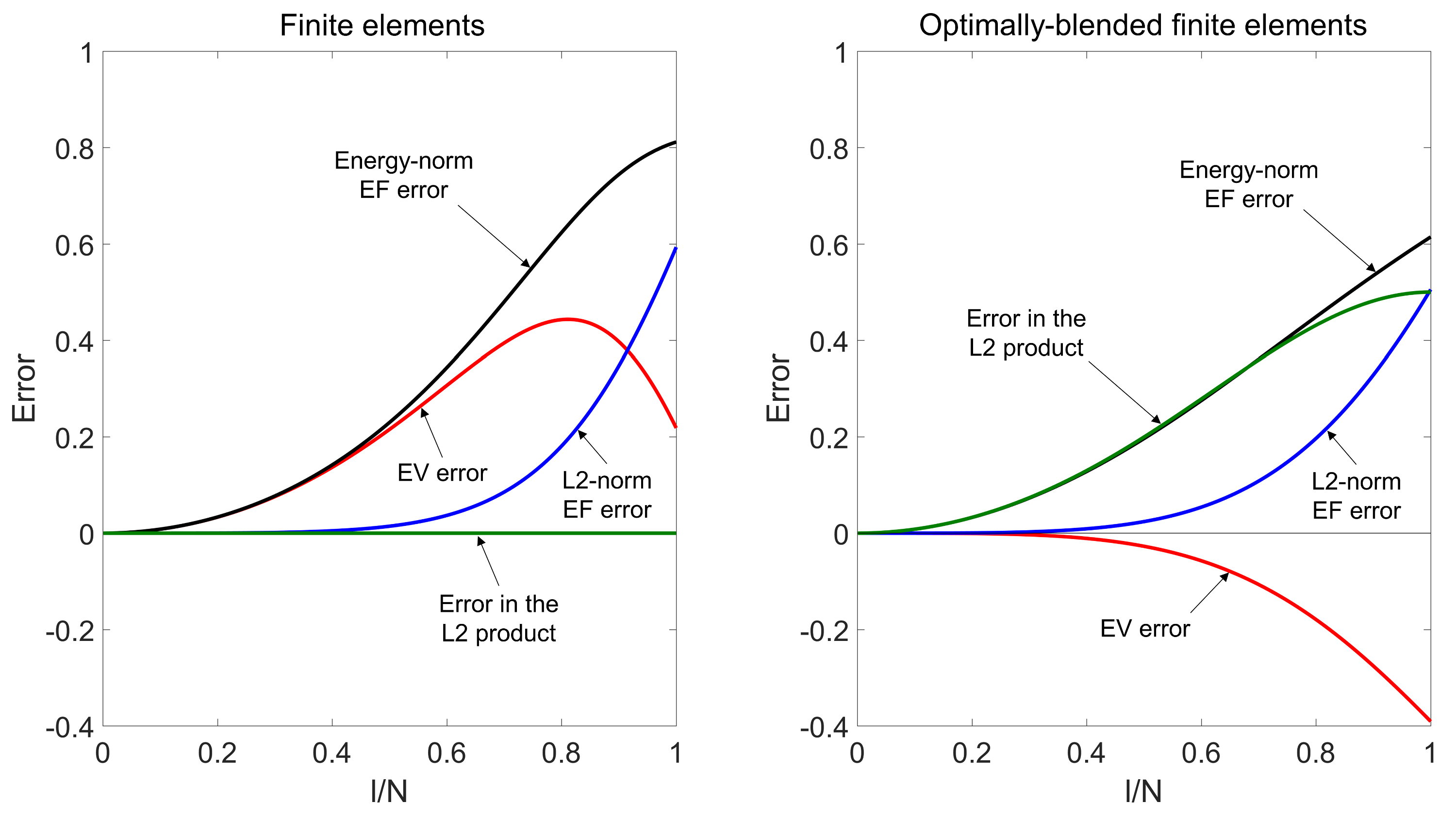}
\caption{Standard linear finite elements (left) versus the optimally-blended spectral-finite element method (right), where EV and EF stand for eigenvalue and eigenfunction, respectively.}
\end{figure}

In the following figures, we show how the use of optimal blending (Section 4) reduces the eigenvalue errors for $C^0$ finite element discretizations. Figures 2 and 3 show the eigenvalue, $L_2$-norm, and energy-norm eigenfunction errors for linear and quadratic $C^0$ elements with the blending parameters $\tau$ equal to $1/2$ and $2/3$, respectively. The fourth term of the modified theorem \eqref{eq:ptfinal} that represents the error in the $L_2$ inner product subtracts from the eigenvalue error thus making this error smaller. The $L_2$- and energy-norm eigenfunction errors are also slightly reduced by the use of optimal quadratures.

\begin{figure}[!ht]
\centering\includegraphics[width=1.0\linewidth]{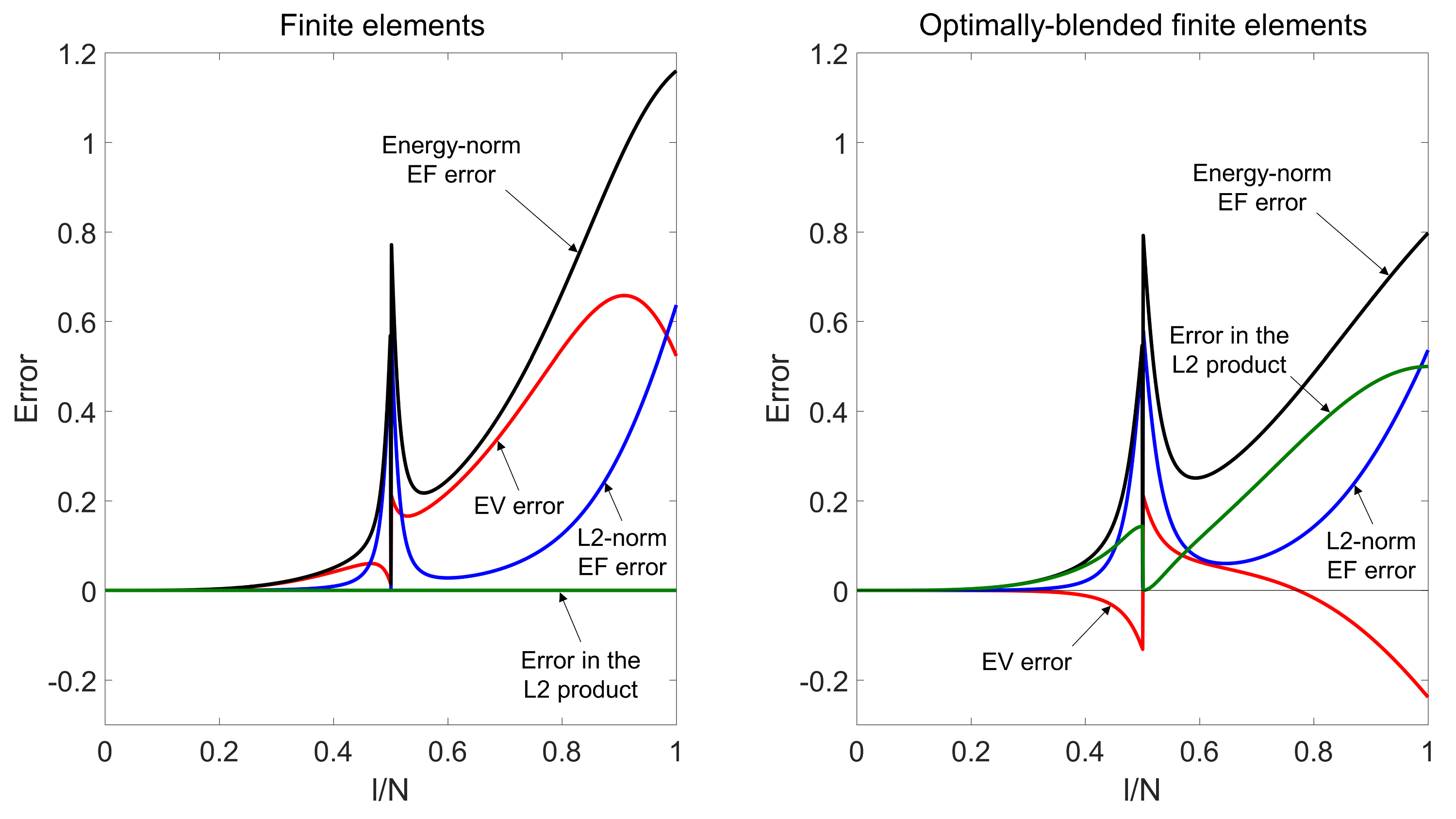}
\caption{Standard quadratic $C^0$ finite elements (left) versus the optimally-blended spectral-finite element method (right), where EV and EF stand for eigenvalue and eigenfunction, respectively.}
\end{figure}

\begin{figure}[!ht]
\centering\includegraphics[width=1.0\linewidth]{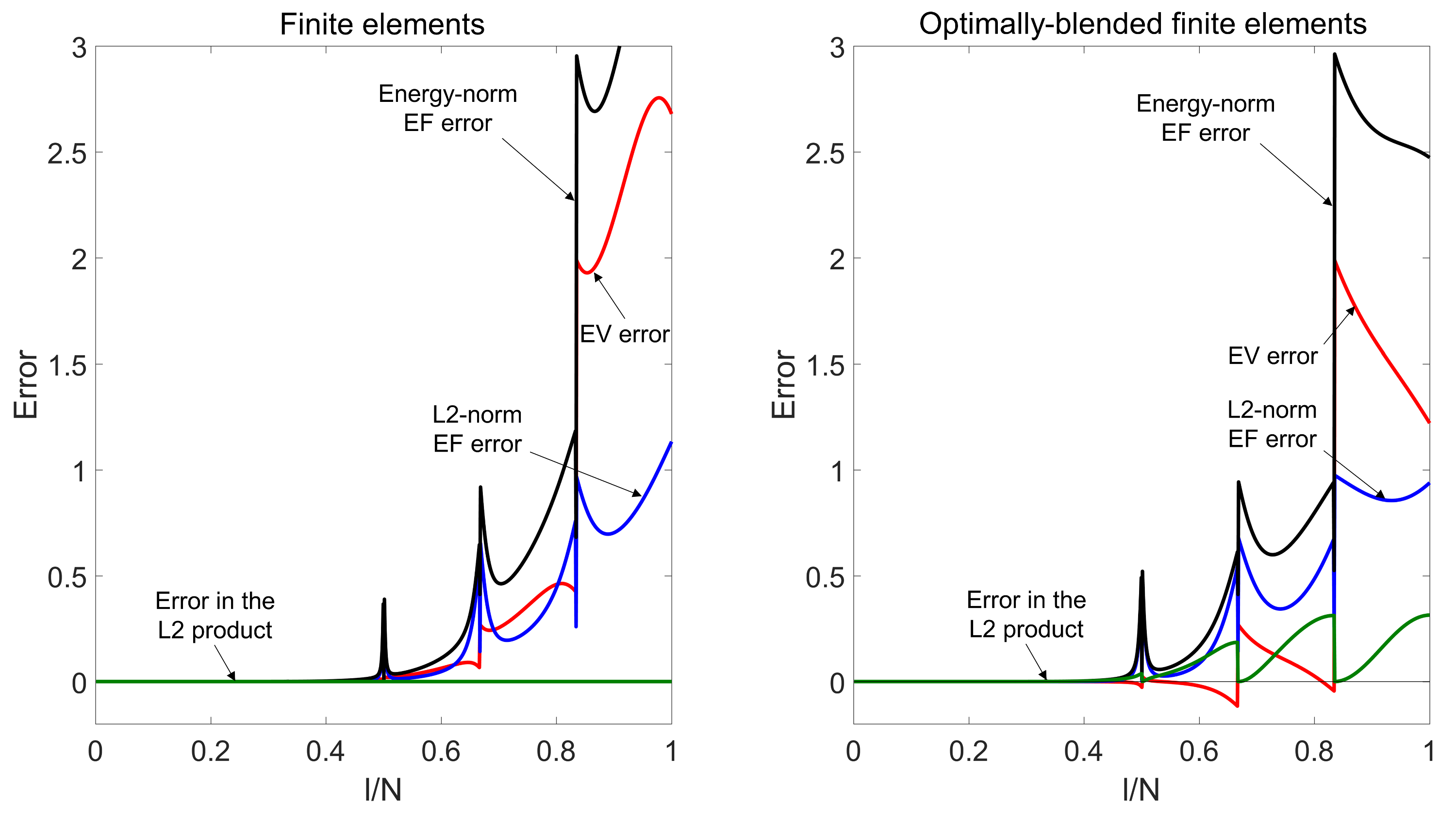}
\caption{Standard sextic $C^0$ finite elements (left) versus the optimally-blended spectral-finite element method (right), where EV and EF stand for eigenvalue and eigenfunction, respectively.}
\end{figure}

The optimal blending improves the spectral properties of higher-order FEA as well. Figure 4 shows the errors for sextic ($p=6, \tau=6/7$) blended finite elements. The y-scale is significantly different in this figure from the previous ones to accommodate the large errors in the optical branch.

\begin{figure}[!ht]
\centering\includegraphics[width=1.0\linewidth]{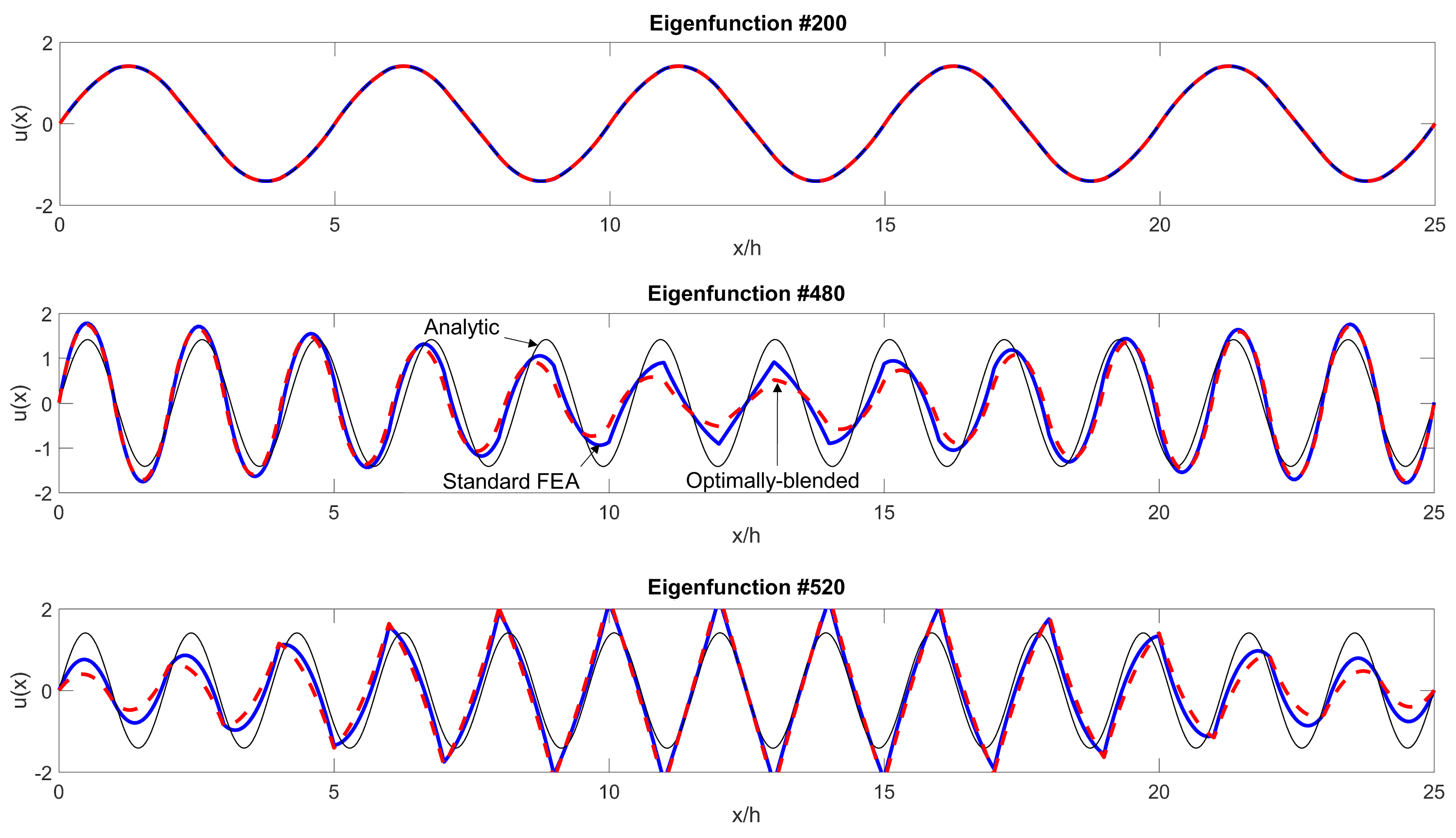}
\caption{Discrete 200th (top), 480th (middle) and 520th (bottom) eigenfunctions for quadratic $C^0$ finite elements (the total number of discrete modes is 1000). The discrete eigenfunctions resulting from the optimally-blended (red) and the standard finite elements (blue) are compared with the analytical eigenfunctions (black).}
\end{figure}

To study the behavior of discrete eigenfunctions from different branches of the spectrum, in Figure 5 we compare the discrete and analytical eigenfunctions for quadratic $C^0$ finite elements. We chose three eigenfunctions: \#200 from the acoustic branch of the spectrum where both the standard and blended finite element schemes have very low error, \#480 and \#520 located near the spike of the stopping band at the center of the spectrum where the errors in the eigenfunctions slightly differ for these methods. As expected, for the lowest mode, both methods provide the results that are qualitatively and quantitatively accurate. On the other hand, the standard finite element method provides slightly better approximation than the optimally-blended scheme for the eigenfunctions located near the stopping-band spike.

\begin{figure}[!ht]
\centering\includegraphics[width=1.0\linewidth]{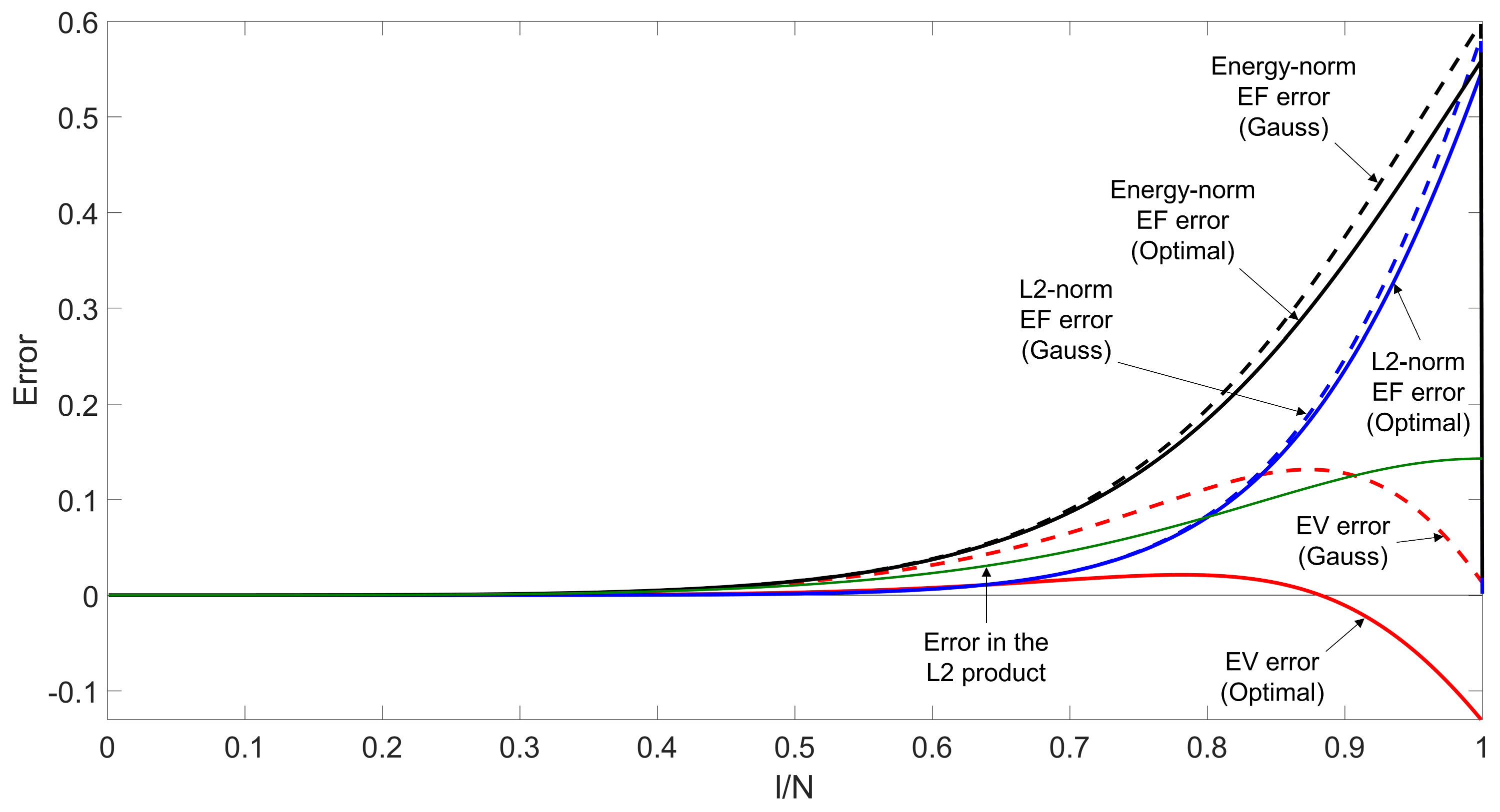}
\caption{Illustration of the Pythagorean eigenvalue theorem for the optimally-blended $C^1$ quadratic elements. EV and EF stand for eigenvalue and eigenfunction, respectively.}
\end{figure}

\begin{figure}[!ht]
\centering\includegraphics[width=1.0\linewidth]{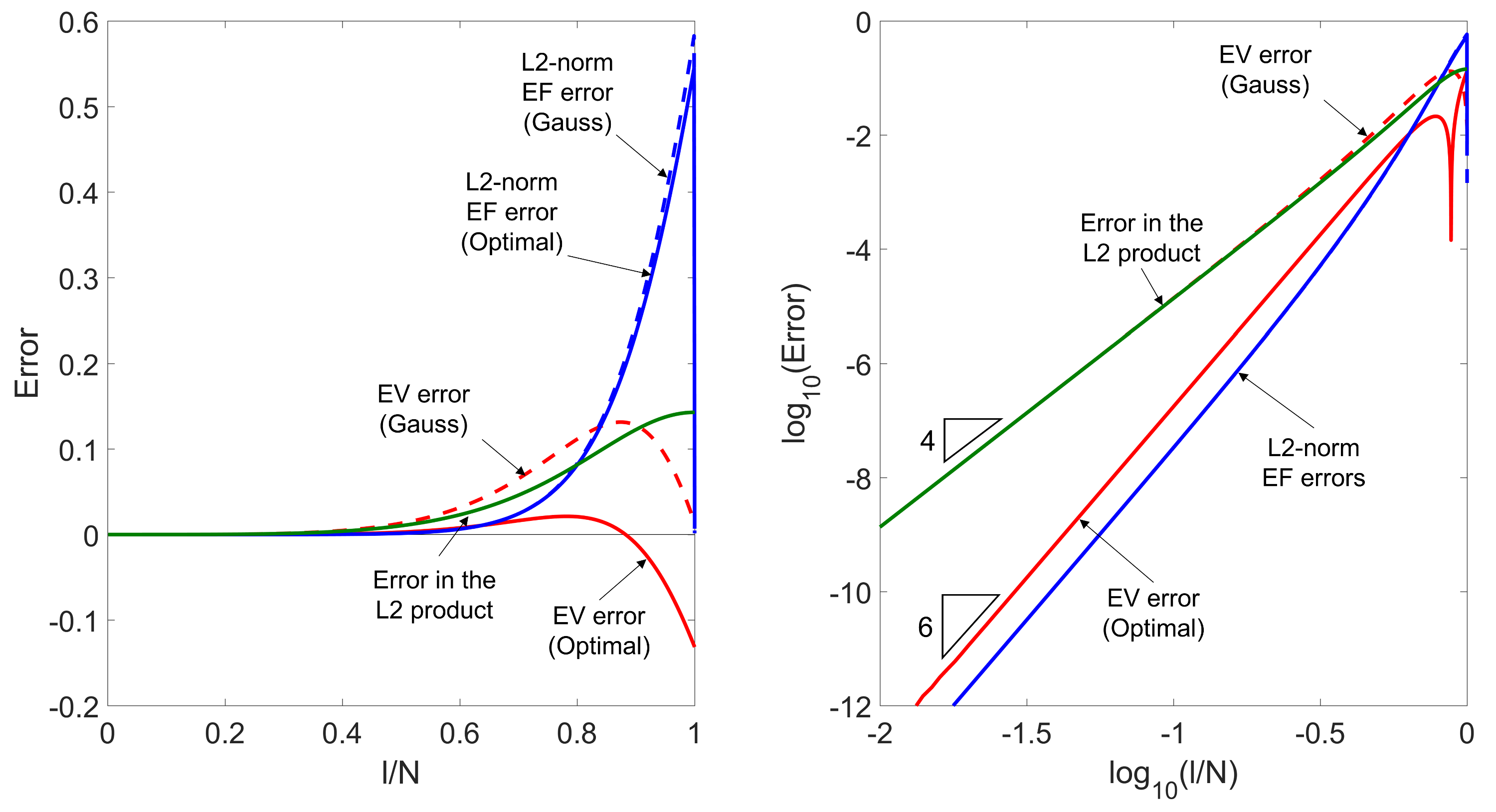}
\caption{Optimal blending for $C^1$ quadratic IGA on linear (left) and logarithmic (right) scales, where EV and EF stand for eigenvalue and eigenfunction, respectively. On the logarithmic scale, the absolute value of the errors is plotted.}
\end{figure}

Let us now show how a blended quadrature rule (a mix of Gauss and Lobatto quadratures in this case) can greatly reduce the eigenvalue errors of isogeometric approximations. Figure 6 shows the eigenvalue, $L_2$-norm, and energy-norm eigenfunction errors for $C^1$ quadratic elements. As can be seen from this figure, the generalized Pythagorean eigenvalue theorem is valid for isogeometric discretizations. The optimal ratio of the Lobatto and Gauss blending is 2 : 1 ($\tau = 2/3$) which is the same ratio proposed by \citet{ainsworth2010optimally} for the FEA case. Figure 7 shows the errors for $C^1$ quadratic elements in both linear and logarithmic scales. With this blending, the scheme has two additional orders of accuracy compared with the standard isogeometric elements. This can be seen in the logarithmic representation of the absolute value of the errors: for the lowest modes (left part of the spectrum), the error in the $L_2$ inner product constitutes the largest part of the total error thus making the eigenvalue error much smaller. The spike in the logarithmic plot is induced by the change in sign of the error. For conciseness, we omit the energy-norm eigenfunction errors from this and the following figures. 

\begin{figure}[!ht]
\centering\includegraphics[width=1.0\linewidth]{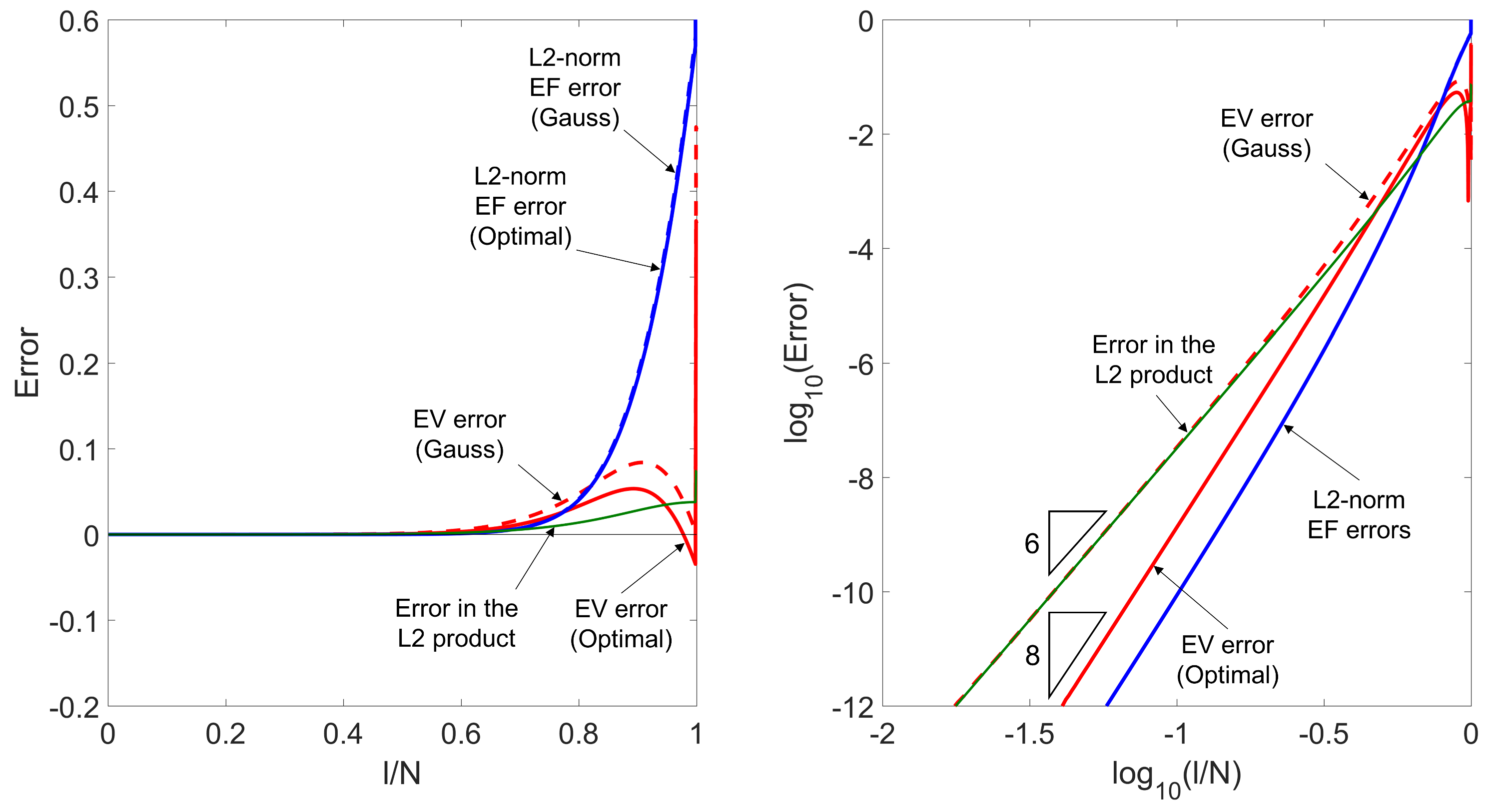}
\caption{Optimal blending for $C^2$ cubic IGA on linear (left) and logarithmic (right) scales, where EV and EF stand for eigenvalue and eigenfunction, respectively. On the logarithmic scale, the absolute value of the errors is plotted.}
\end{figure}

Figure 8 shows the errors in the optimally-blended $C^2$ cubic IGA scheme. The optimal in this case is a non-convex blending of $-3/2$ Gauss and $5/2$ Lobatto rules \citep{deng2017quadratures}. We observe much smaller eigenvalue errors compared to the standard $C^2$ isogeometric elements. The convergence rate of the optimally-blended IGA scheme is $O(\Omega^8)$. The figure shows the numerical evidence and \citep{deng2017quadratures} shows the analysis of this blending.

\begin{figure}[!ht]
\centering\includegraphics[width=1.0\linewidth]{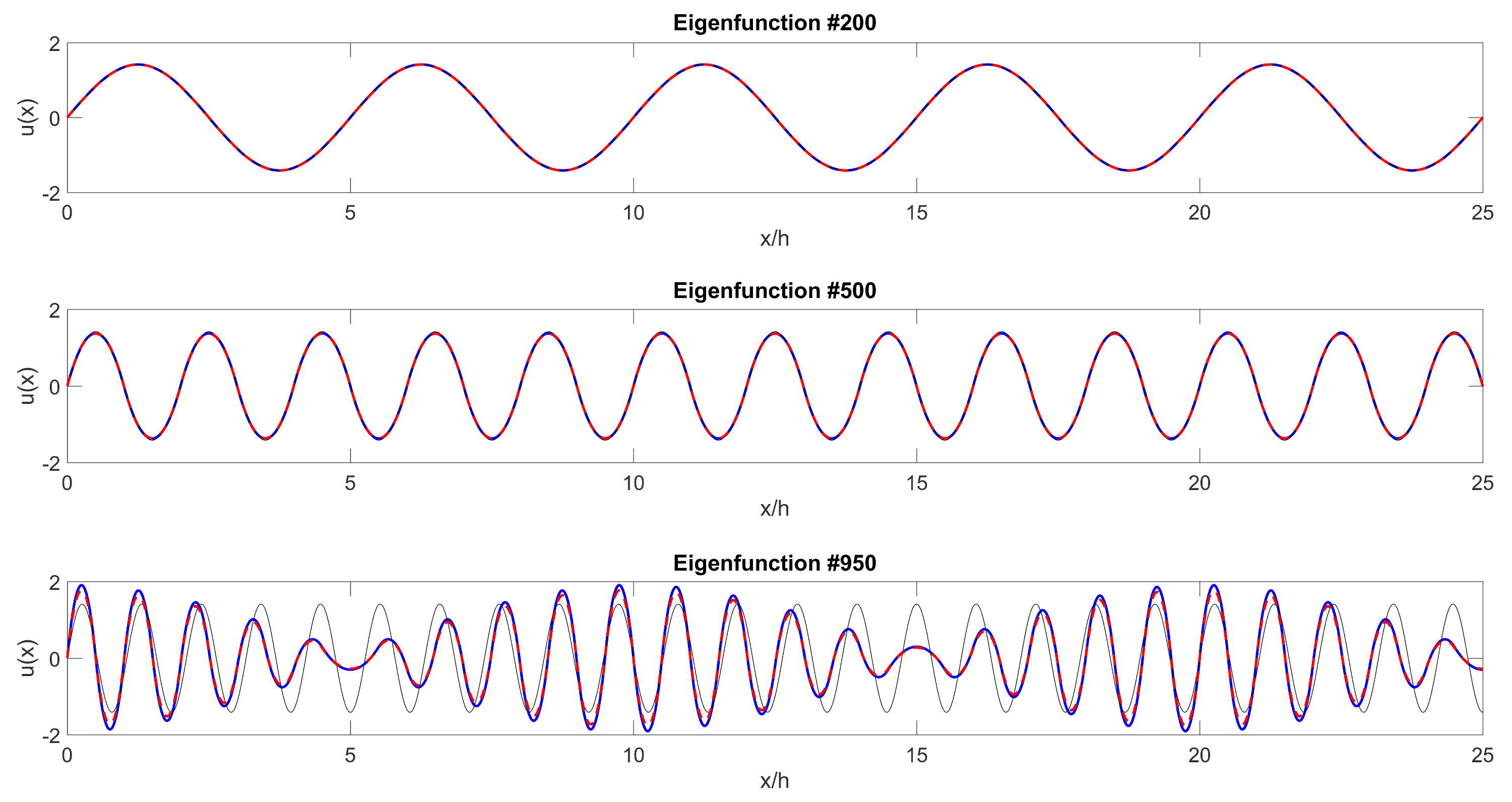}
\caption{Discrete 200th (top), 500th (middle) and 950th (bottom) eigenfunctions for $C^1$ quadratic elements (the total number of discrete modes is 1000). The discrete eigenfunctions resulting from the optimally-blended (red) and the standard isogeometric elements (blue) are compared with the analytical eigenfunctions (black).}
\end{figure}

These results demonstrate that the blended smooth isogeometric inner products significantly improve the accuracy of the discrete approximations when compared to the fully-integrated Gaussian counterparts. The fourth term of the modified Pythagorean eigenvalue theorem \eqref{eq:ptfinal} represents the "correction" to the approximate eigenvalues. The $L_2$- and energy-norm eigenfunction errors are almost not affected by the blending. The logarithmic format confirms the theoretical convergence is achieved for blended $C^2$ elements and the generalized Pythagorean theorem is valid to high precision in the numerical experiments.

The eigenfunction error, contrary to the finite element case, is either slightly improved (for low-order IGA) or almost not affected at all (for high-order elements). Three discrete and analytical eigenfunctions for $C^1$ quadratic elements are shown in Figure 9. Again, for the 200th mode, both blended and non-blended schemes are very accurate. The approximation of the higher 500th eigenfunction is also good, while the 800th is less accurate. Both the standard and blended schemes provide similar results and no loss of accuracy in eigenfunction approximation is observed due to the use of the non-standard blended quadratures.

\begin{figure}[!ht]
\centering\includegraphics[width=1.0\linewidth]{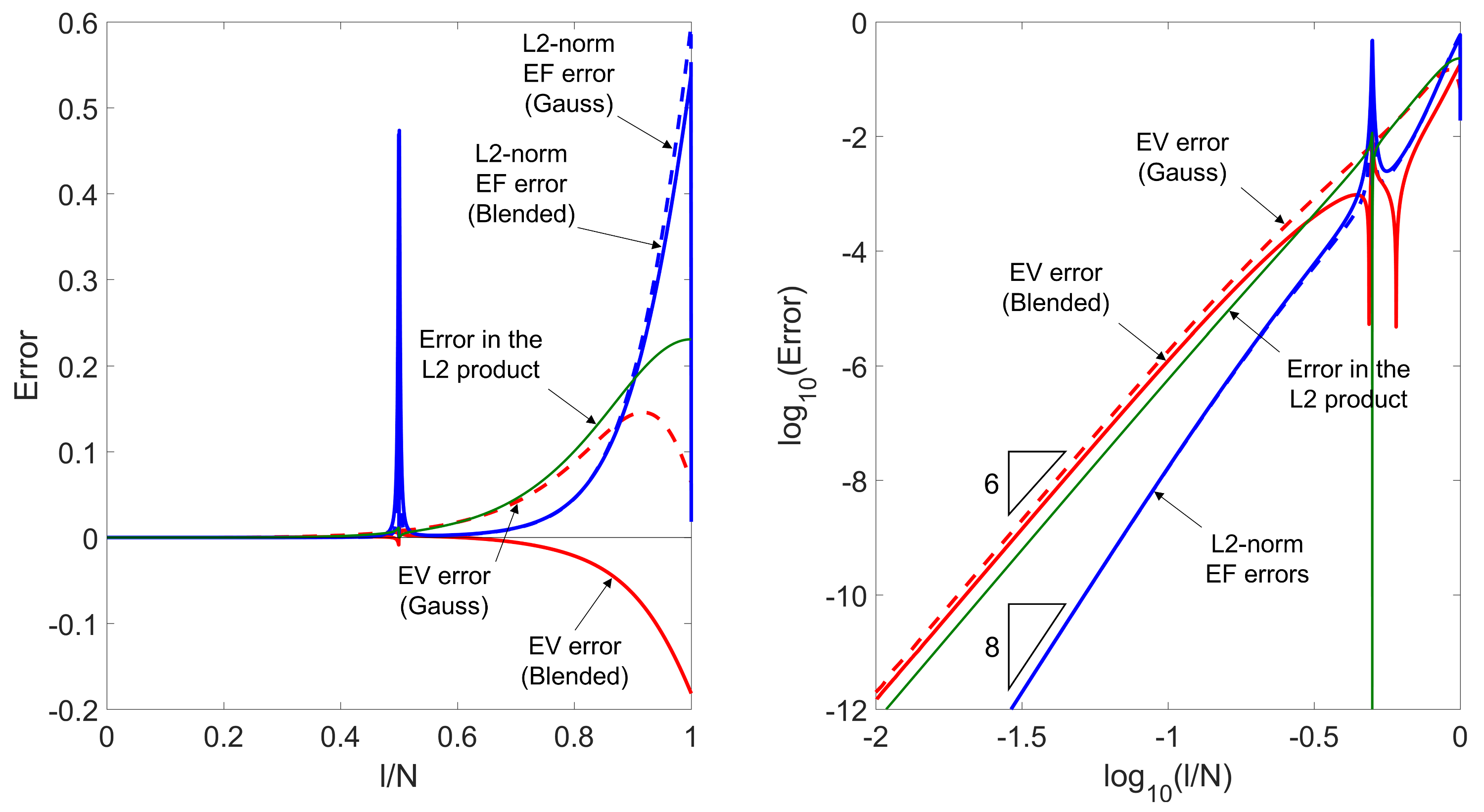}
\caption{Blending with $\tau=3/4$ for $C^1$ cubic elements on linear (left) and logarithmic (right) scales, where EV and EF stand for eigenvalue and eigenfunction, respectively. On the logarithmic scale, the absolute value of the errors is plotted.}
\end{figure}

Now, we consider the spectral error analysis for blended isogeometric methods with less continuous spaces (less than $p-1$ continuity). Figure 10 compares the blending with $\tau=3/4$ for $C^1$ cubic elements versus the standard fully integrated scheme. This blending would be optimal for $C^0$ elements of the same order, however, in this case it does not improve the convergence rate of the method. Nevertheless, a positive impact on the eigenvalue errors can be observed in the whole spectrum, especially in its middle part.

\begin{figure}[!ht]
\centering\includegraphics[width=1.0\linewidth]{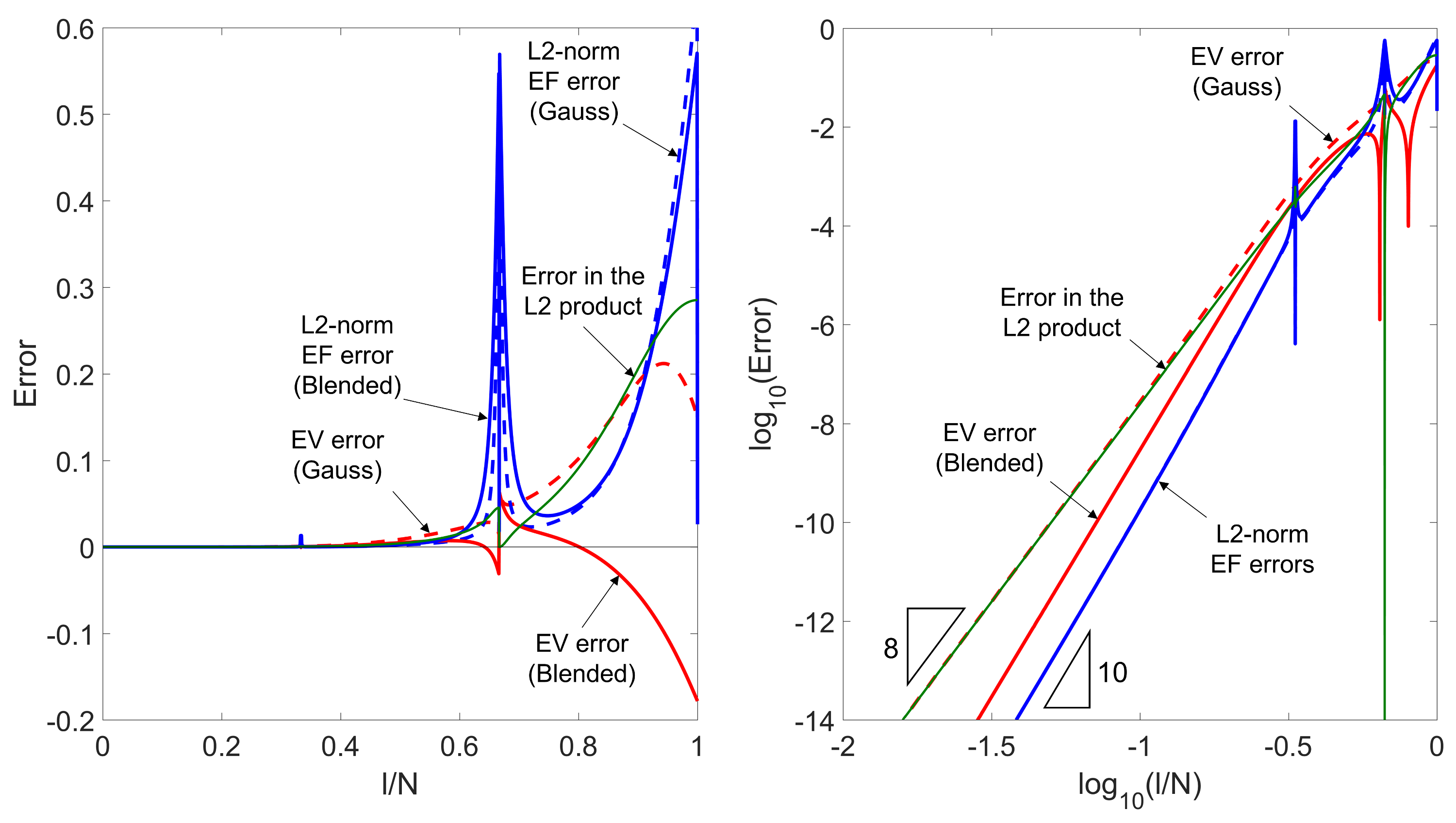}
\caption{Blending with $\tau=4/5$ for $C^1$ quartic elements on linear (left) and logarithmic (right) scales, where EV and EF stand for eigenvalue and eigenfunction, respectively. On the logarithmic scale, the absolute value of the errors is plotted.}
\end{figure}

Errors for $C^1$ quartic elements are shown in Figure 11. The $C^0$ optimal blending scheme ($\tau=4/5$) used in this case has a convergence order of 10, i.e. this blending is able to deliver two additional orders of convergence compared to 8 of the standard quartic elements. Thus, this scheme is optimal in the Ainsworth's sense \cite{ainsworth2010optimally}. The eigenvalue errors are significantly reduced in the whole spectrum.

\begin{figure}[!ht]
\centering\includegraphics[width=1.0\linewidth]{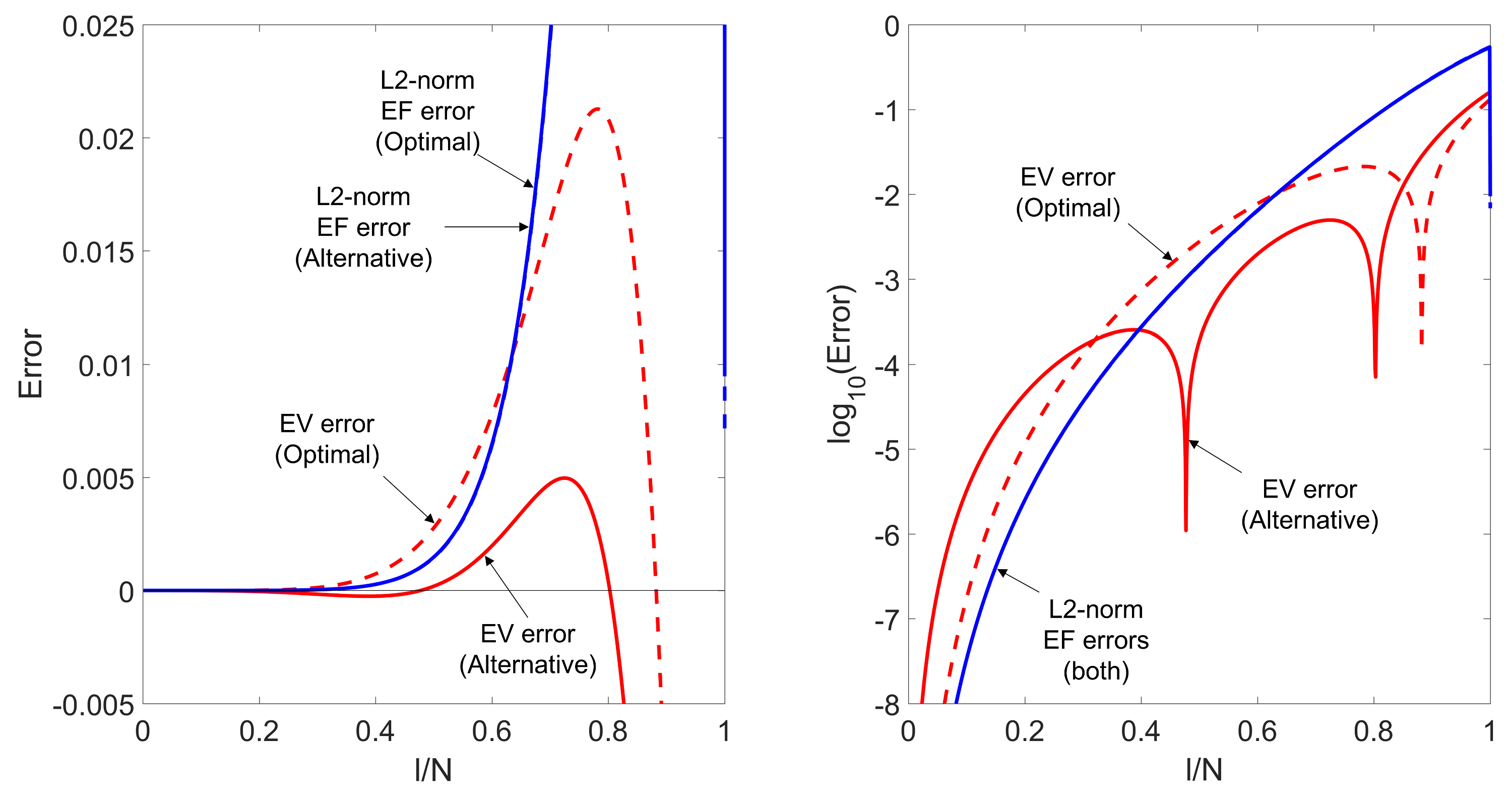}
\caption{Alternative $\tau=5/6$ blending versus the optimal $\tau=2/3$ scheme for $C^1$ quadratic IGA on linear (left) and linear-logarithmic (right) scales, where EV and EF stand for eigenvalue and eigenfunction, respectively.}
\end{figure}

\subsection{Alternative blending schemes}

In the following examples, we show the errors for quadratic and cubic isogeometric elements using alternative blendings and compare them with the optimal ones. As stated in Section 5, for practical applications, we may seek acceptable errors for desired intervals of wavenumber (frequency) for a fixed mesh size. Figure 12 shows the effect on the spectral distribution of the error for the blending parameter $\tau=5/6$ applied on a $C^1$ quadratic discretization. While this blending is not optimal, in the sense that it does not deliver superconvergence, it results in eigenvalue errors that are much smaller than the "optimal" ones in the range $[0.34\ 0.96]$.

\begin{figure}[!ht]
\centering\includegraphics[width=1.0\linewidth]{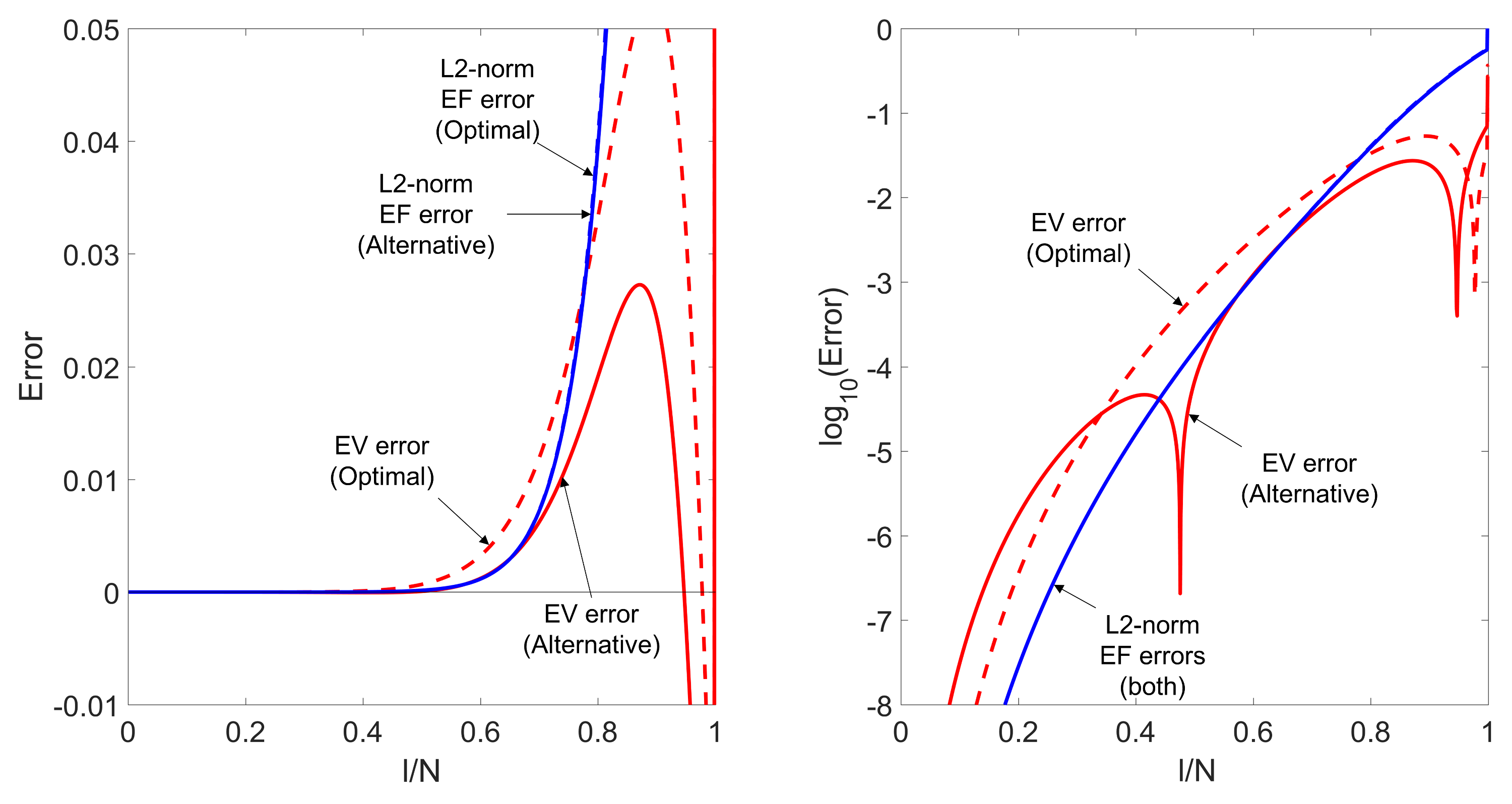}
\caption{Alternative $\tau=5$ blending versus the optimal $\tau=5/2$ scheme for $C^2$ cubic IGA on linear (left) and linear-logarithmic (right) scales, where EV and EF stand for eigenvalue and eigenfunction, respectively.}
\end{figure}

Figure 13 compares the blending with $\tau=5$ for cubic $C^2$ elements versus the optimal blending. Again, the optimal blending has much smaller errors for the lowest modes, but the alternative blending has better approximation properties at the higher ones that are of practical interest in wave propagation problems.

\begin{figure}[!ht]
\centering\includegraphics[width=1.0\linewidth]{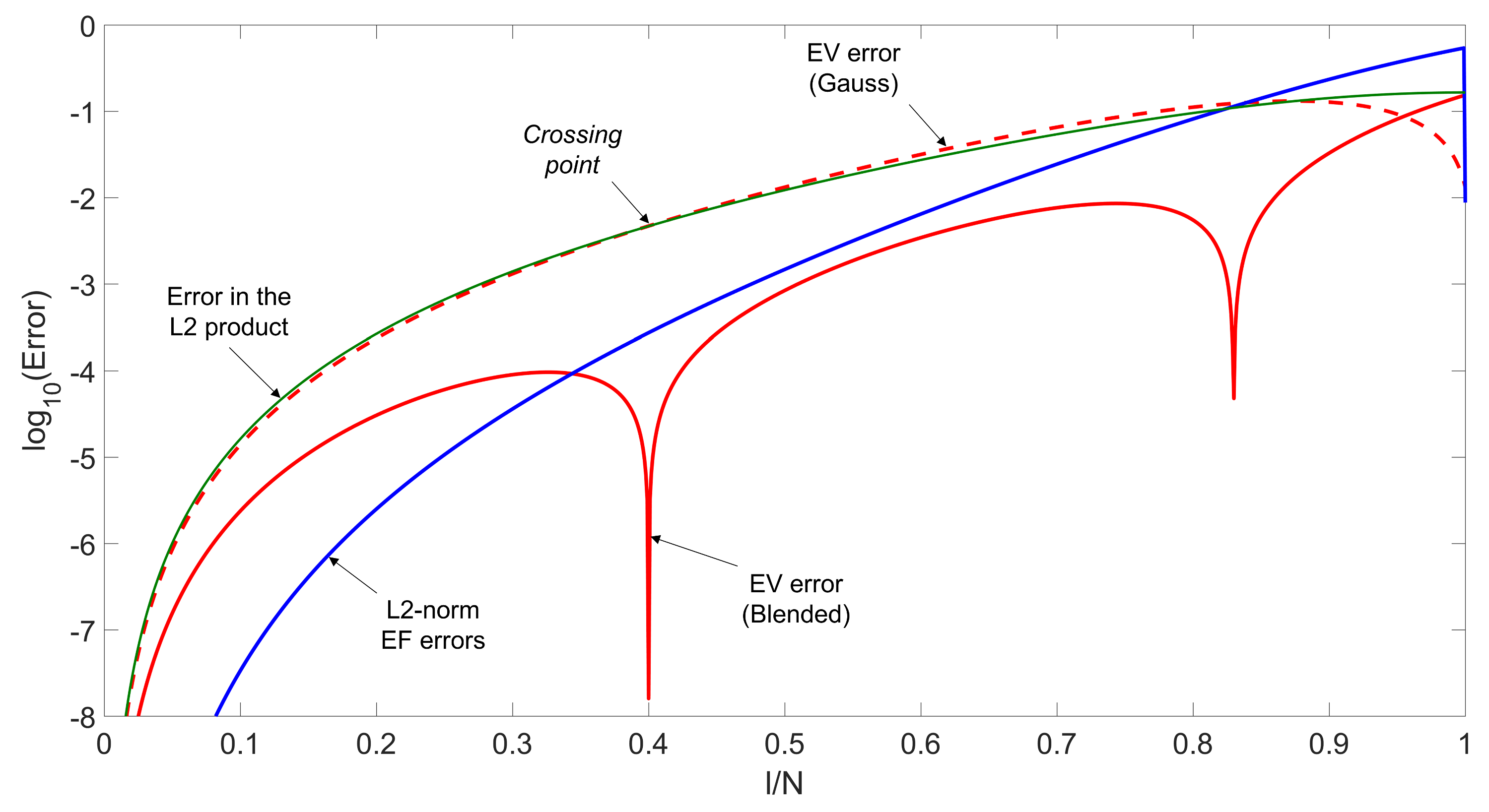}
\caption{$C^1$ quadratic elements blended with $\tau=0.20895$. The error at $l/N = 0.4$ is reduced almost to zero, where EV and EF stand for eigenvalue and eigenfunction, respectively.}
\end{figure}

The blending parameter can be chosen to deliver a better approximation for a particular, user-specified value of the wavenumber. Assume that we wish the error to vanish at $\frac{k^h h}{\pi} = 0.4$ (five points per wavelength). By using the blending parameter $\tau=0.20895$, we obtain a numerical scheme based of $C^1$ quadratic elements whose approximation properties are shown in Figure 14.

\begin{figure}[!ht]
\centering\includegraphics[width=1.0\linewidth]{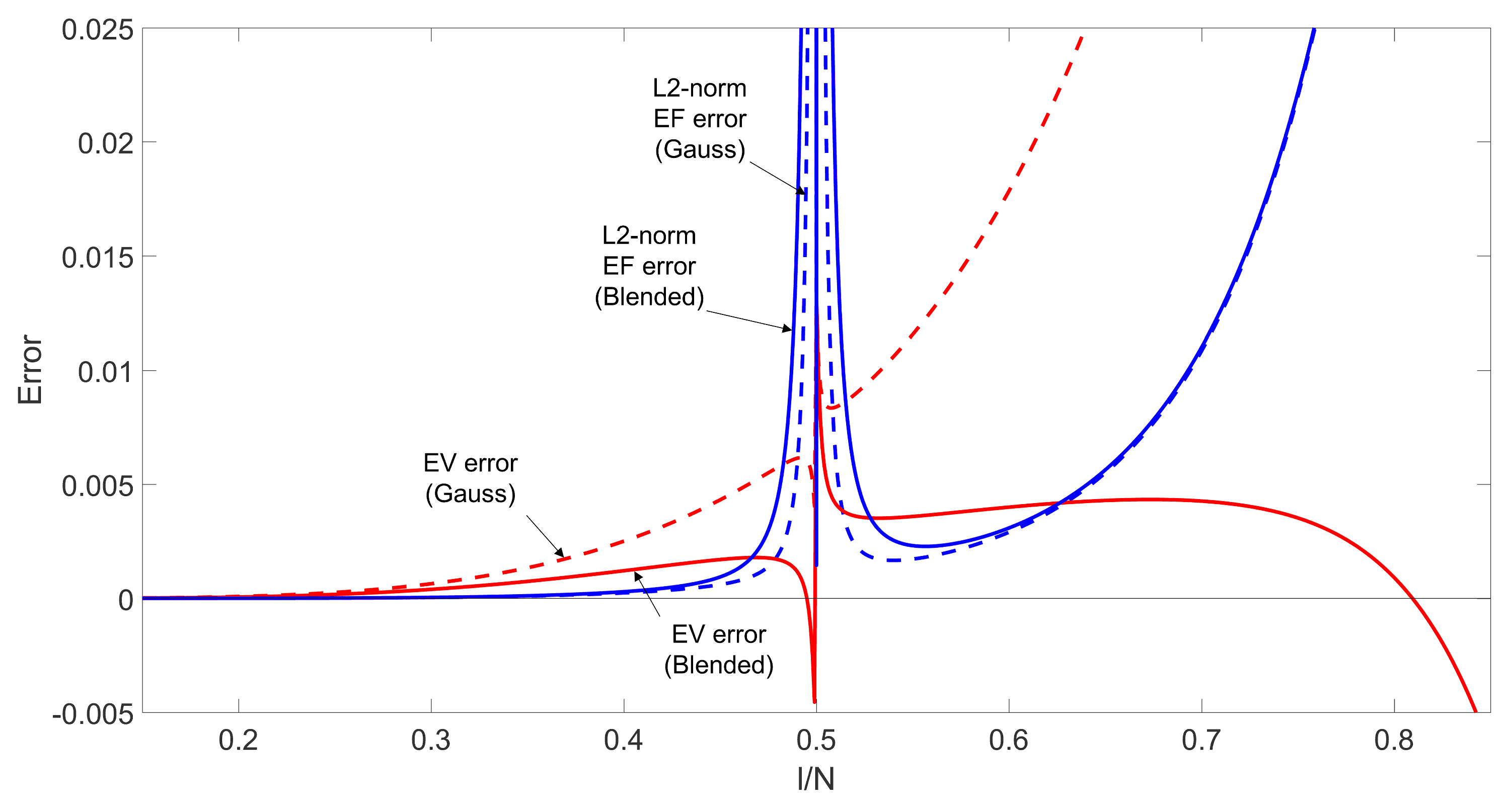}
\caption{The middle part of the spectrum for $C^1$ cubic elements blended with $\tau=7/12$, where EV and EF stand for eigenvalue and eigenfunction, respectively.}
\end{figure}

Alternative blendings can be also successfully used for isogeometric elements with reduced continuity. Figure 15 shows the blending scheme with $\tau=7/12$ for $C^1$ cubic elements. This blending greatly reduces the eigenvalue errors in the middle range; for example, the errors in the first 80\% of the spectrum are within 2\% range.

Summarizing the results above, the optimally-blended $C^0$ and $C^{p-1}$ elements exhibit much smaller eigenvalue errors and have two additional orders of superconvergence when compared with the standard schemes based on a Gauss quadrature rule that fully integrates all the inner products. Moreover, a properly chosen non-optimal blending rule can reduce the eigenvalue errors even further for a range of normalized eigenvalues (or set the error to zero for a particular mode) while preserving the convergence rate of the standard discretizations.

\subsection{Optimally-blended methods in 2D}

Now we continue our study with the dispersion properties of the two-dimensional eigenvalue problem / wave equation on tensor product meshes. An optimal blending parameter for a multidimensional problem coincides with the optimal parameter for the one dimensional case \citep{ainsworth2010optimally}. The exact eigenvalues and eigenfunctions of the 2D eigenvalue problem are given by
\begin{equation}
{\lambda _{jk}} = {(j^2+k^2)}{\pi ^2},\ \ \ {u_{jk}} = 2 \sin (j\pi x) \sin (k\pi y),
\end{equation}
for $j,k = \overline {1,\infty }$. Again, the approximate eigenvalues $\lambda _{jk}^h$ are sorted in ascending order. The resulting 1D plots of eigenvalue errors show oscillations which are common in this kind of analysis \citep{reali2004isogeometric, hughes2008duality}.

\begin{figure}[!ht]
\centering\includegraphics[width=1.0\linewidth]{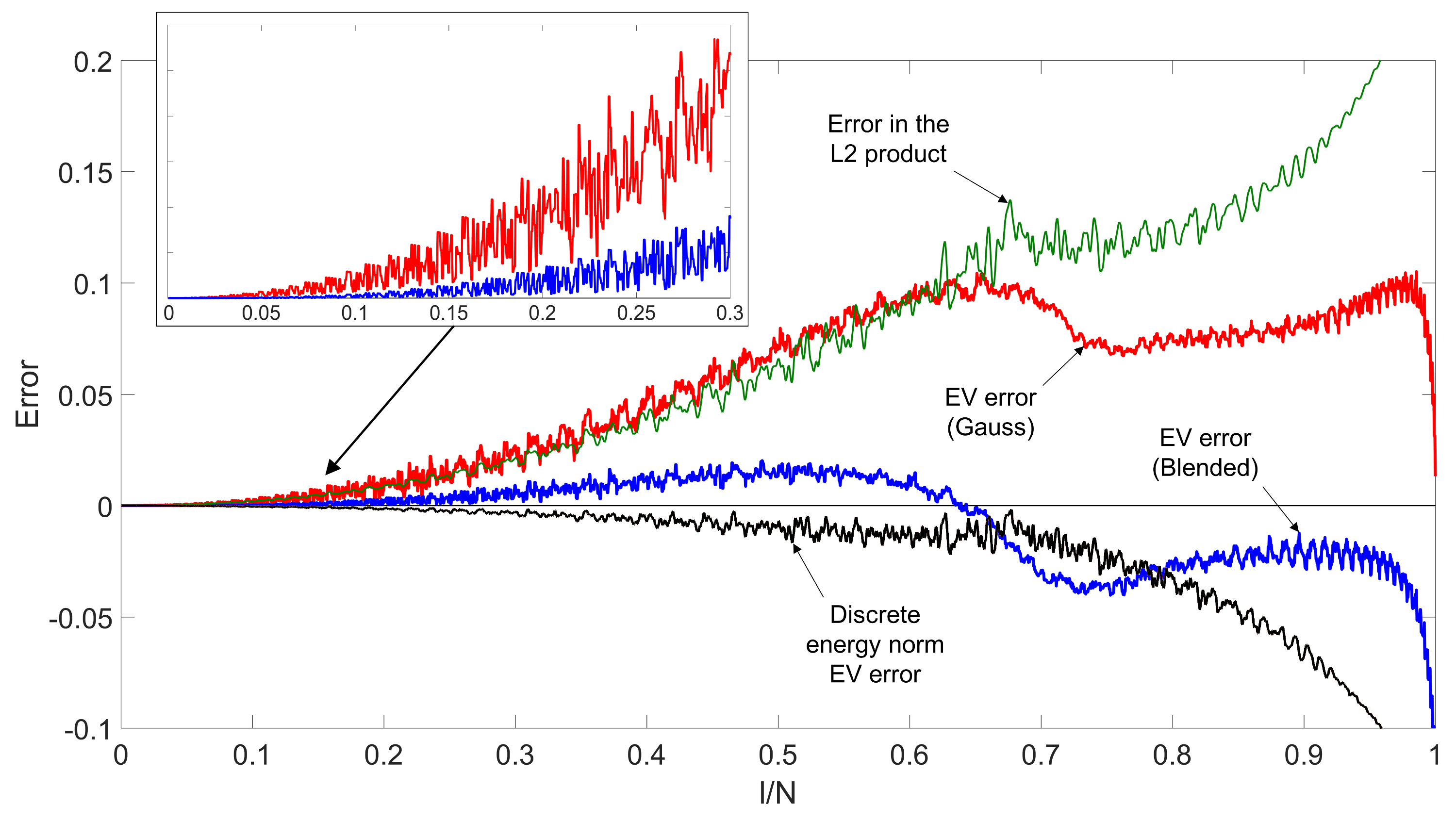}
\caption{Eigenvalue errors for $C^1$ quadratic elements using the standard and optimally-blended schemes. The inset compares the accuracy in the first 30\% of the spectrum. EV and EF stand for eigenvalue and eigenfunction, respectively. Now the eigenvalue errors for the Gauss and blended scheme are shown in red and blue colors, respectively.}
\end{figure}

Figure 16 compares the eigenvalue errors in the standard 1D format for the $C^1$ quadratic elements using the standard Gauss and the optimally-blended scheme. Now both the third and the fourth terms of the modified Pythagorean eigenvalue theorem are not zero and modify the eigenvalue approximation. The  error in the discrete energy norm is detrimental to the approximation when it has a negative sign. Some of its detrimental contribution is offset by the error in the $L_2$ norm. Nevertheless, the optimally-blended scheme has significantly smaller approximation errors in the whole spectrum.

\begin{figure}[!ht]
\centering\includegraphics[width=1.0\linewidth]{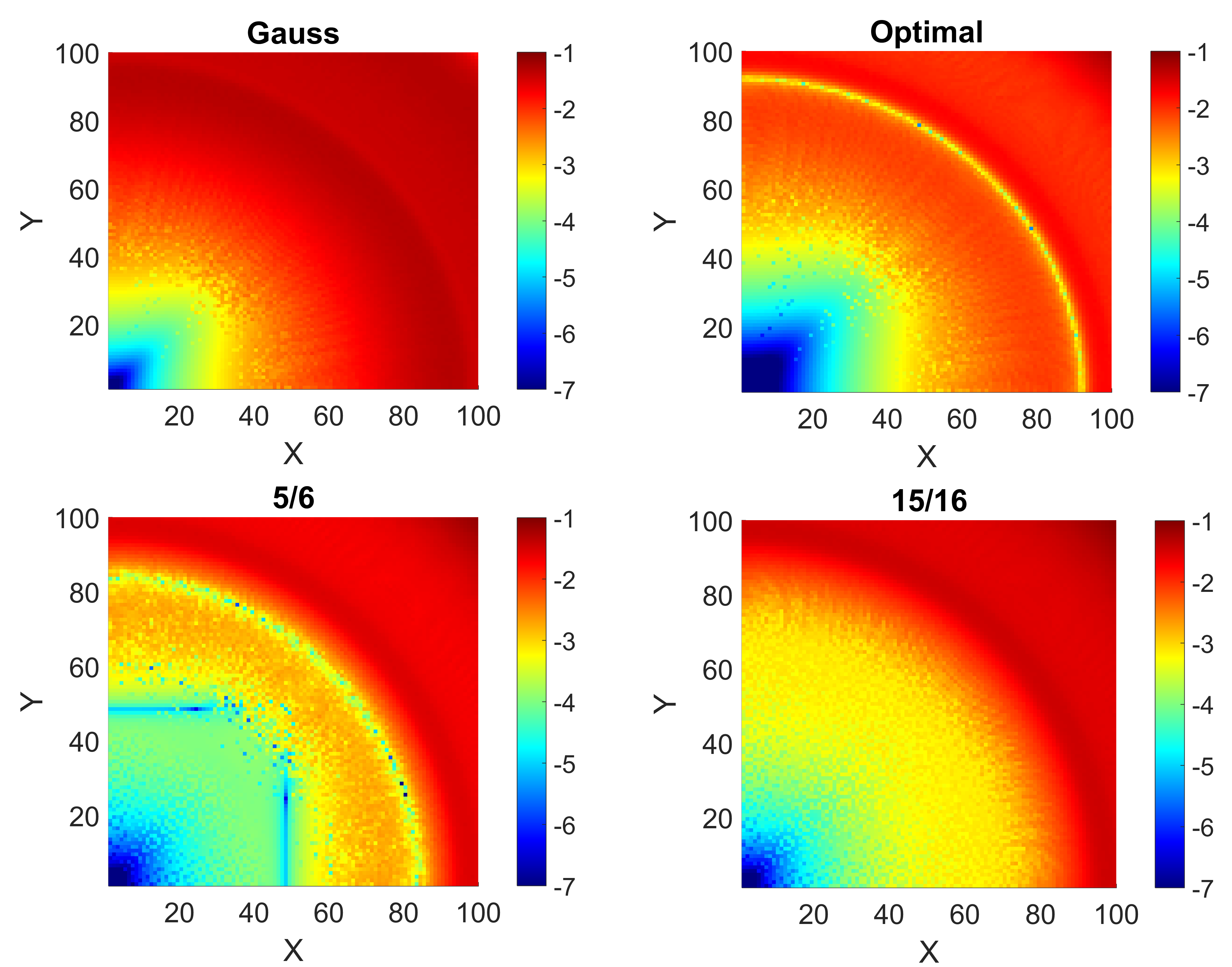}
\caption{Different blendings for $C^1$ quadratic elements. Color represents the logarithm of the absolute value of the eigenvalue error.}
\end{figure}

Figure 17 compares the eigenvalue errors of the fully integrated by Gauss inner products for $C^1$ quadratic elements with the optimally-blended scheme ($\tau=2/3$) and two alternative blendings ($\tau=5/6$ and $\tau=15/16$) in a 2D format. We can observe that for 2D problems the main features of the optimal and alternative non-optimal blendings are preserved: the optimal one has smaller errors for the lowest modes, but the alternative ones have better approximation properties in the middle part of spectra.

\begin{figure}[!ht]
\centering\includegraphics[width=1.0\linewidth]{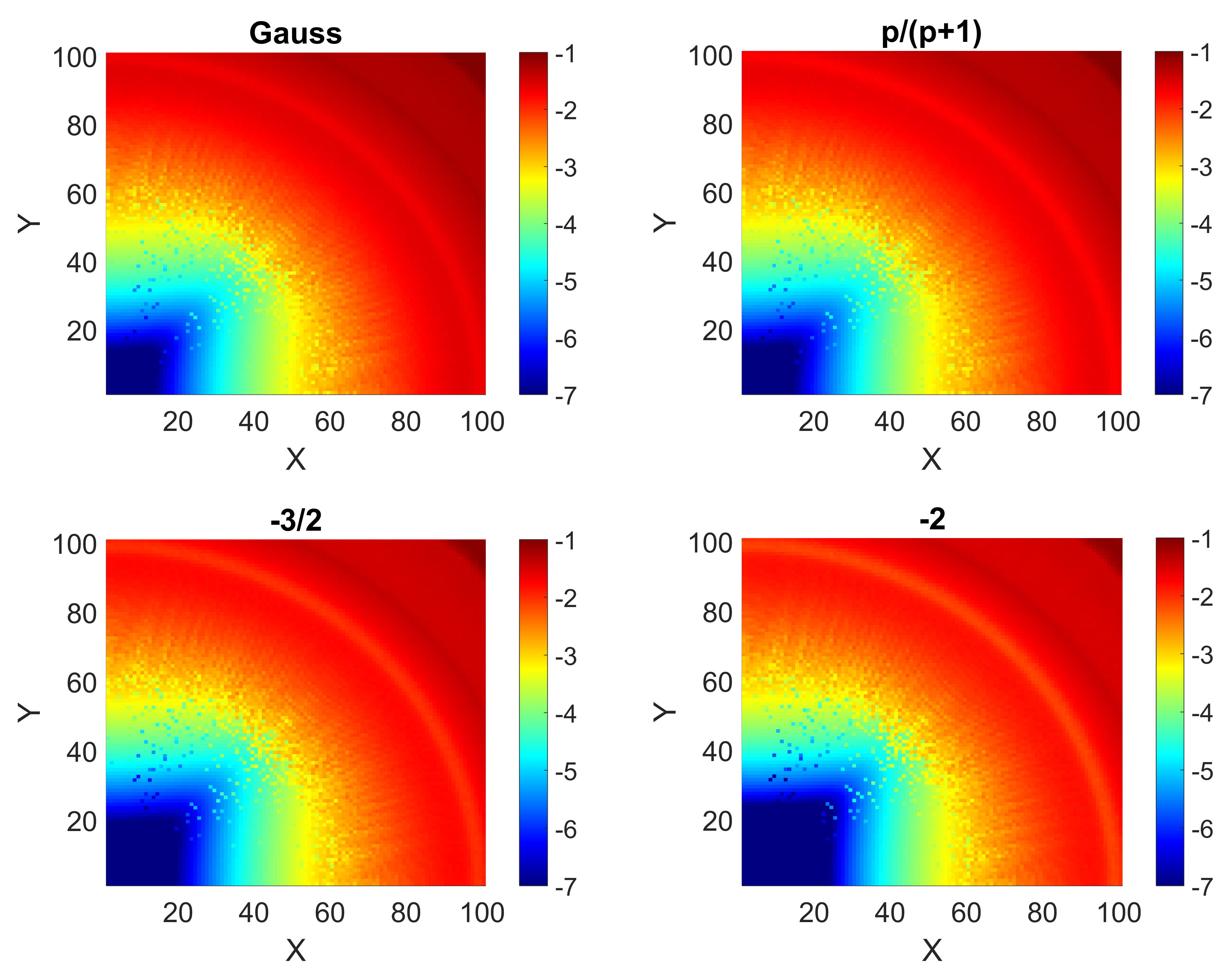}
\caption{Different blendings for $C^2$ cubic elements. Color represents the logarithm of the absolute value of the eigenvalue error.}
\end{figure}

Figure 18 shows the impact on the simulation results for a $C^2$ cubic basis when a $p+1$ Gauss quadrature is compared against several blending schemes (FEA optimal value of $\tau=3/4$, IGA optimal value of $\tau=5/2$, and an alternative blending with $\tau=3$). As can be seen from these results, using a value of the blending parameter $\tau=3/4$, which is optimal for $C^0$ cubic elements, does not change the approximation properties compared to the standard scheme. On the contrary, the non-convex blendings lead to smaller eigenvalue errors.

\subsection{Non-uniform grids}

Until now, we have considered only uniform grids that are commonly used in dispersion analysis. However, real-life applications involve modeling problems that include elements of very different shapes and sizes. In this section, we study the efficiency of the blended quadratures on non-uniform grids.

\begin{figure}[!ht]
\centering\includegraphics[width=1.0\linewidth]{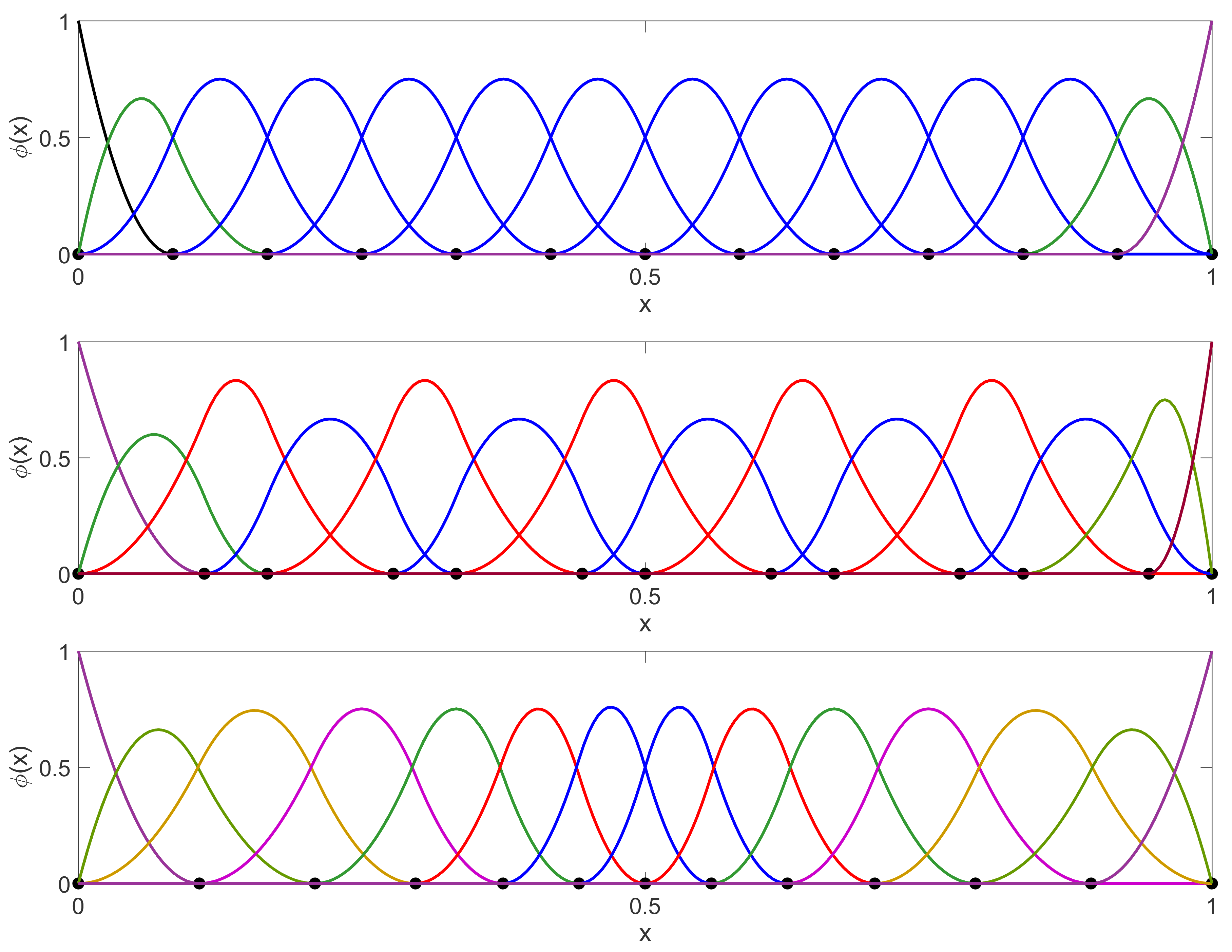}
\caption{Example of a $C^1$ quadratic basis for a uniform grid with 12 elements (top), a grid with two types of elements (middle), and a stretched grid (bottom)}
\end{figure}

Using non-uniform knot spacing leads to different types of basis function to be present in the system. Figure 19 shows the $C^1$ quadratic basis for three different grid refinements (boundary knots are shown as one). In the first case, the uniform grid refinement (the size of all elements is $h$) leads to a homogeneous basis everywhere, except at the boundaries of the domain. In the second case, a mesh that combines elements of two sizes, $\frac{2}{3}h$ and $\frac{4}{3}h$ (where $h$ is the mesh size of an equivalent uniform mesh with the same number of elements), leads to two families of basis functions. In the third case, we consider a classical stretched grid that is refined such that
\begin{equation}
{h}_{k+1} = \alpha {h}_{k},
\end{equation}
and the mesh is refined towards the center of the domain. The basis on this grid has a large variety of different functions.

\begin{figure}[!ht]
\centering\includegraphics[width=1.0\linewidth]{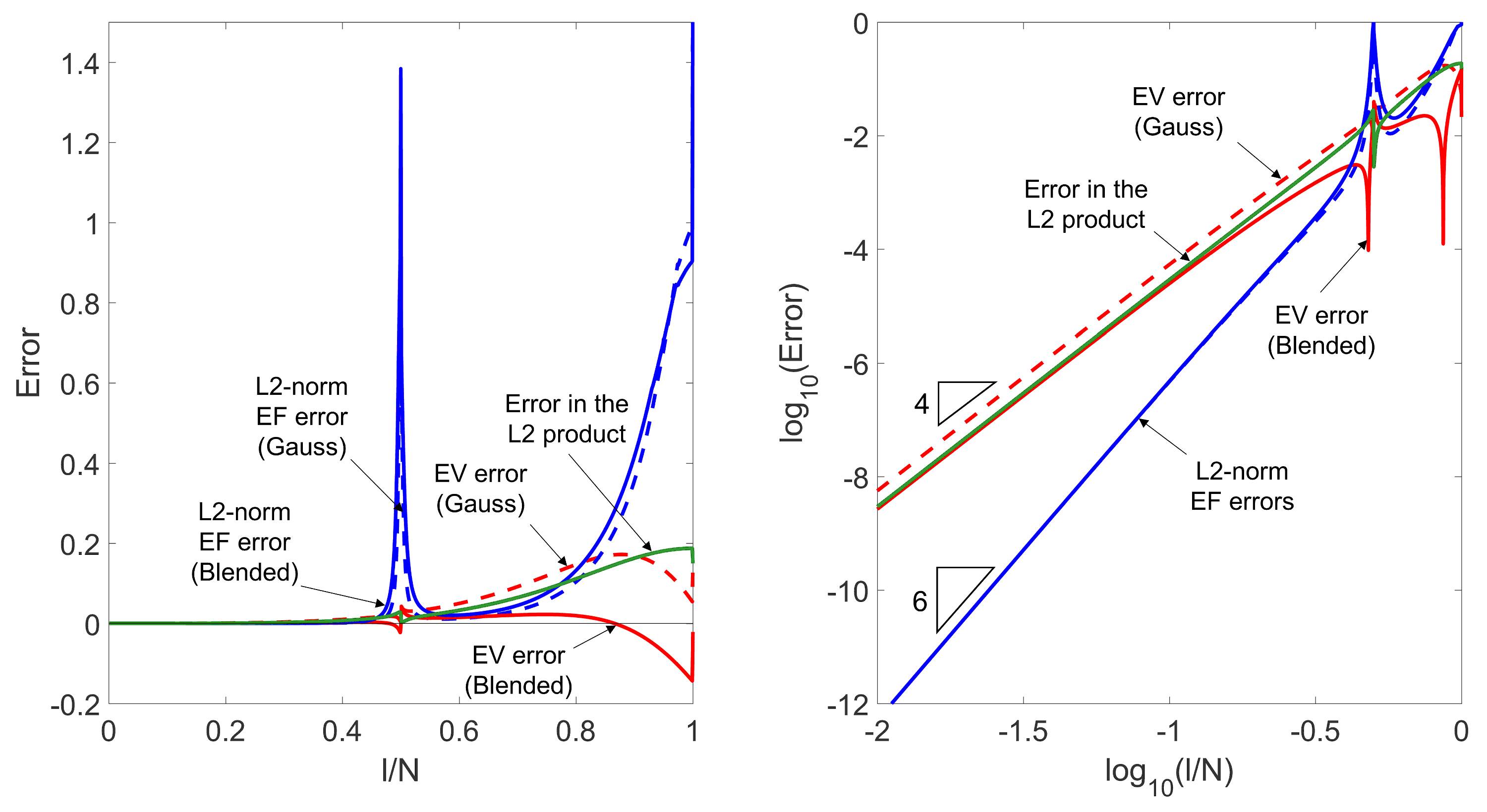}
\caption{Gauss versus optimal blending using $C^1$ quadratic elements on the grid with $\frac{2}{3}h$ and $\frac{4}{3}h$ elements, where EV and EF stand for eigenvalue and eigenfunction, respectively.}
\end{figure}

Figure 20 compares the Gauss quadrature with the optimal blending for $C^1$ quadratic basis on a uniform grid on the grid that combines elements of two sizes, $\frac{2}{3}h$ and $\frac{4}{3}h$. The main feature of the spectrum is the large spike in the eigenfunction error at the center. The inner basis functions of the problem belong to two different groups and this, similar to the finite element case, leads to branching of the spectrum (despite the fact that the functions have the same continuity). The use of the blended quadrature leads to smaller eigenvalue errors thus confirming the efficiency of the method. The convergence is optimal as can be seen from the logarithmic representation of the errors, but no superconvergence is observed when the mesh is non-uniform. Nevertheless, we can determine numerically another quadrature that leads to two additional orders of superconvergence for this mesh. Figure 21 shows the errors for a non-convex blending of three-point Gauss and Lobatto quadratures with $\tau \approx 1.27$.

\begin{figure}[!ht]
\centering\includegraphics[width=1.0\linewidth]{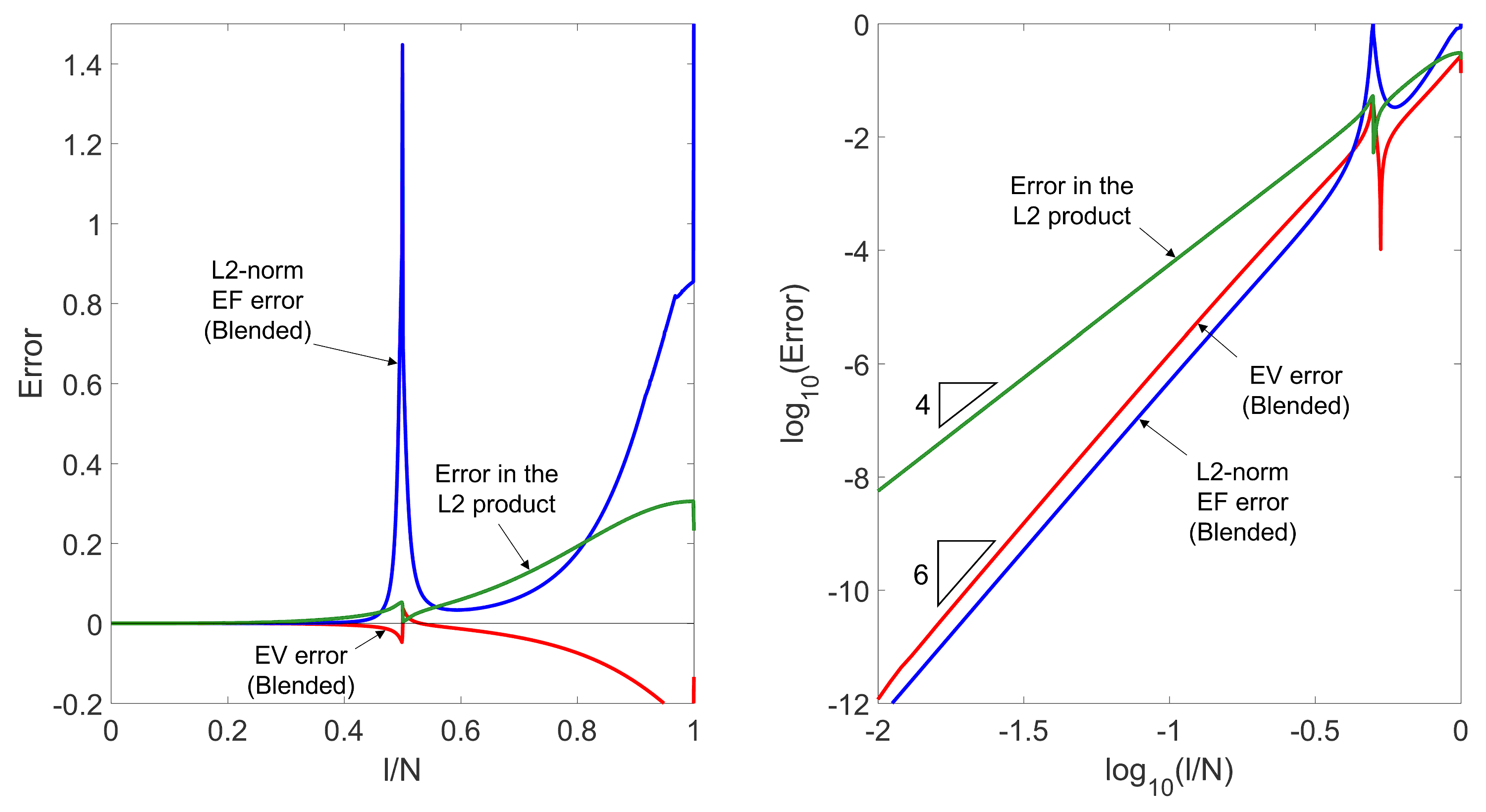}
\caption{Alternative blending with $\tau \approx 1.27$ for $C^1$ quadratic elements on the grid with $\frac{2}{3}h$ and $\frac{4}{3}h$ elements, where EV and EF stand for eigenvalue and eigenfunction, respectively.}
\end{figure}

Finally, Figure 22 shows the results for a stretched grid with the stretching parameter $\alpha$ equal to 1.02. Again, the use of the optimal $\tau=\frac{p}{p+1}$ quadrature leads to much more accurate results. Surprisingly, not only the eigenvalues, but also the eigenfunctions of the problem are better approximated in this case.

\begin{figure}[!ht]
\centering\includegraphics[width=0.7\linewidth]{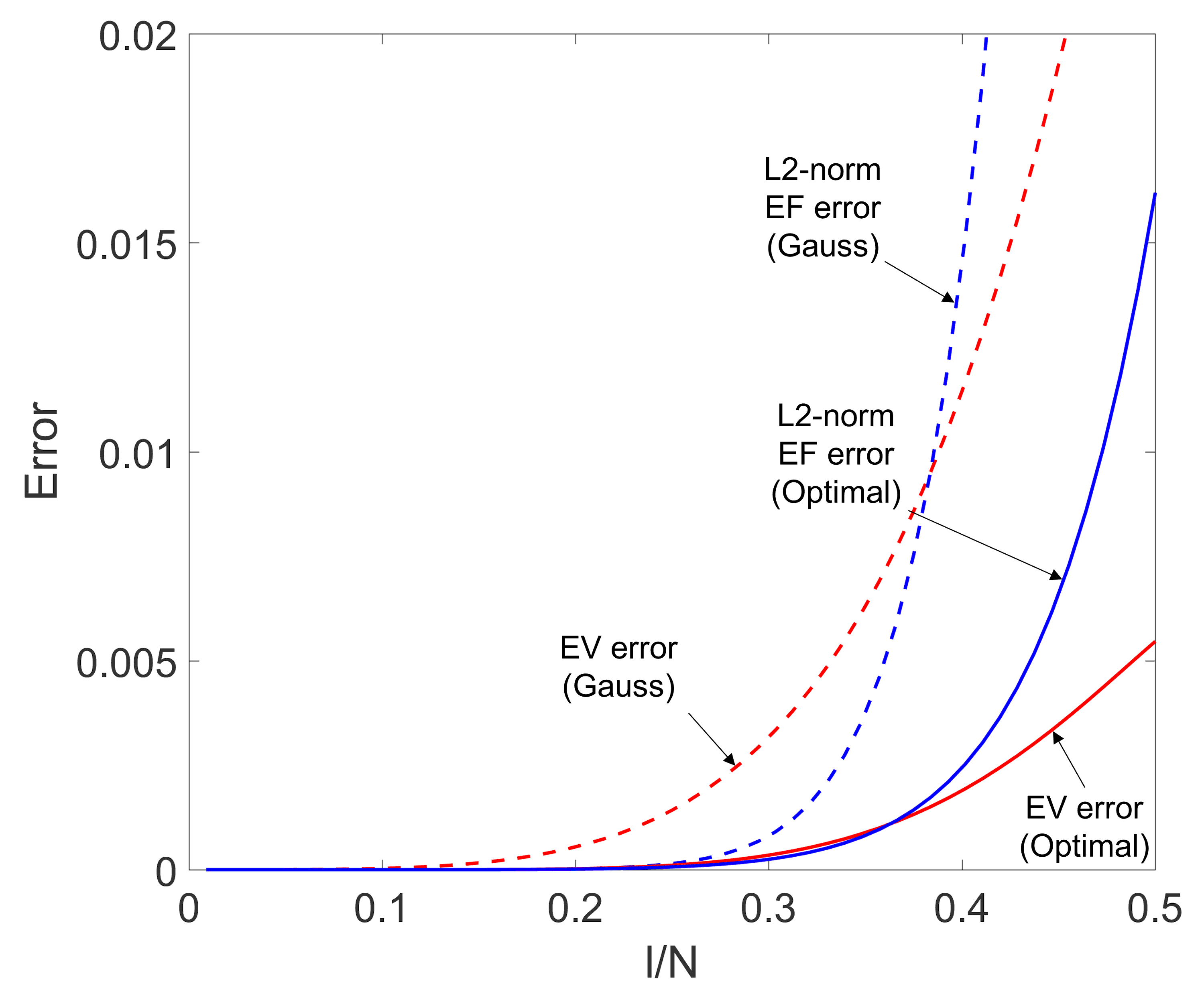}
\caption{Gauss versus optimal blending for $C^1$ quadratic elements on the stretched grid with $\alpha = 1.02$, where EV and EF stand for eigenvalue and eigenfunction, respectively.}
\end{figure}





\section{Conclusions and future outlook}

We show that using blended quadrature rules reduces the phase error of the numerical method, without affecting the overall efficacy of the method. To explain the observed behavior, we generalize the Pythagorean eigenvalue theorem to account for the effects of the modified inner products on the resulting weak forms. The proposed technique improves the superior spectral accuracy of isogeometric analysis. We can extend the method to arbitrary high-order $C^{p-1}$ isogeometric elements by identifying suitable quadrature rules. So far, equivalent quadrature rules, which are not the result of blending a Gauss and a Lobatto quadrature, for arbitrary isogeometric elements are not available in the existing literature and will be the subject of our future work.

Another future direction of research will focus on isogeometric discretizations with variable continuity. We will study how the breaks in continuity and inhomogeneity of the basis they produce affect the dispersion properties of the method and how these effects can be minimized by the use of blended quadratures.

\section{Acknowledgments}

This publication was made possible in part by a National Priorities Research Program grant 7-1482-1-278 from the Qatar National Research Fund (a member of The Qatar Foundation), and by the European Union's Horizon 2020 Research and Innovation Program of the Marie Sklodowska-Curie grant agreement No. 644202. The J. Tinsley Oden Faculty Fellowship Research Program at the Institute for Computational Engineering and Sciences (ICES) of the University of Texas at Austin has partially supported the visits of VMC to ICES. The Spring 2016 Trimester on ``Numerical methods for PDEs'', organized  with the collaboration of the Centre Emile Borel at the Institut Henri Poincare in Paris supported VMC's visit to IHP in October, 2016.

\section*{References}

\nocite{*}

\bibliographystyle{elsarticle-harv}\biboptions{square,sort,comma,numbers}
\bibliography{references}








\end{document}